\documentclass[reqno]{amsart}
\usepackage{amssymb}
\usepackage{amsmath}

\makeatletter
\@addtoreset{equation}{section}
\makeatother

\renewcommand\thefigure{\thesection.\@arabic\c@figure}
\renewcommand\thetable{\thesection.\@arabic\c@table}

\newtheorem{theorem}{Theorem}[section]
\newtheorem{lemma}[theorem]{Lemma}
\newtheorem{proposition}[theorem]{Proposition}
\newtheorem{corollary}[theorem]{Corollary}

\newcommand{\mc}[1]{{\mathcal #1}}
\newcommand{\mf}[1]{{\mathfrak #1}}
\newcommand{\mb}[1]{{\mathbf #1}}
\newcommand{\bb}[1]{{\mathbb #1}}
\newcommand{\bs}[1]{{\boldsymbol #1}}

\newcommand{\<}{\langle}
\renewcommand{\>}{\rangle}

\begin{document}

\title[Incompressible Navier-Stokes equation]
{A lattice gas model for the incompressible
  Navier-Stokes equation}

\author{J. Beltr\'an, C.~Landim}

\address{\noindent IMCA, Calle los Bi\'ologos 245, Urb. San C\'esar
  Primera Etapa, Lima 12, Per\'u and PUCP, Av. Universitaria cdra. 18,
  San Miguel, Ap. 1761, Lima 100, Per\'u. 
\newline e-mail: \rm
  \texttt{johel@impa.br} }

\address{\noindent IMPA, Estrada Dona Castorina 110, CEP 22460 Rio de
  Janeiro, Brasil and CNRS UMR 6085, Universit\'e de Rouen, UMR 6085,
  Avenue de l'Universit\'e, BP.12, Technop\^ole du Madrillet, F76801
  Saint-\'Etienne-du-Rouvray, France.
\newline e-mail: \rm
  \texttt{landim@impa.br} }

\begin{abstract}
  We recover the Navier-Stokes equation as the incompressible limit of
  a stochastic lattice gas in which particles are allowed to jump over
  a mesoscopic scale. The result holds in any dimension assuming the
  existence of a smooth solution of the Navier-Stokes equation in a
  fixed time interval. The proof does not use non-gradient methods or
  the multi-scale analysis due to the long range jumps.
\end{abstract}

\maketitle

\section{introduction}

A major open problem in non-equilibrium statistical mechanics is the
derivation of the hydrodynamical equations from microscopic
Hamiltonian dynamics. The main difficulty in this project lies in the
poor knowledge of the ergodic properties of such systems. To overcome
this obstacle, deterministic Hamiltonian dynamics have been
successfully replaced by interacting particle systems (cf.
\cite{kl} and references therein).

Following this approach, in the sequel of the development of the
non-gradient method by Quastel \cite{q} and Varadhan \cite{v},
Esposito, Marra and Yau \cite{emy, emy2} derived the incompressible
Navier-Stokes equation for stochastic lattice gases in dimension $d\ge
3$. 

The main step of their proof relies on a sharp estimate of the
spectral gap of the jump part of the generator of the process and on
the characterization of the germs of the exact and closed forms in a
Hilbert space of local functions. The characterization of the closed
forms as the sum of exact forms and currents allows, through
a multi-scale analysis, the decomposition of the current as a sum of a
gradient part and a local function in the range of the generator.

In this article we consider a stochastic lattice gas with long range
jumps. The dynamics is build in a way that the density and the
momentum are the only conserved quantities.  Choosing appropriately
the size and the rates of the jumps, we are able to show that a small
perturbation of a constant density and momentum profile evolves in a
diffusive time scale as the solution of the incompressible
Navier-Stokes equation.

In contrast with \cite{emy, emy2}, the mesoscopic range of the jumps
permits to consider perturbations around the constant profile of order
$N^{-b}$, for $b$ small, where $N$ is a scaling parameter proportional
to the inverse of the distance between particles. This choice has two
important consequences. On the one hand, in order to close the
equation, one does not need to replace currents by averages of
conserved quantities over macroscopic boxes, but only over mesoscopic
cubes, whose size depend on the parameter $b$. In particular, there is
no need to recur to the multi-scale analysis or to the closed and
exact forms, simplifying considerably the proof. On the other hand,
choosing $b$ small enough ($b<1/2$), one can avoid in dimension $1$
and $2$ the Gaussian fluctuations around the hydrodynamic limit and
prove a law of large numbers for the conserved quantities in this
regime. We are thus able to derive the incompressible Navier-Stokes
equation even in low dimension, where the usual approach is intrinsically
impossible since it involves scales in which fluctuations appear.

The main drawback of the approach presented is that it requires a
bound on the spectral gap of the full dynamics restricted to finite
cubes. The bound needs only to be polynomial in the volume of the
cube, but the generator includes the collision part. This problem,
already mentioned in \cite{emy2}, is rather difficult in general. We
prove such a bound in Section \ref{gap} for a specific choice of
velocities.

The model can be informally described as follows.  Let $\mc V$ be a
finite set of velocities in $\bb R^d$, invariant under reflections and
exchange of coordinates. For each $v$ in $\mc V$, consider a
long-range asymmetric exclusion process on $\bb Z^d$ whose mean drift
is $v N^{-(1+b)}$. Superposed to this
dynamics, there is a collision process which exchange velocities of
particles in the same site in a way that momentum is conserved.

Under diffusive time scaling, assuming local equilibrium, it is not
difficult to show that the evolution of the conserved quantities is
described by the parabolic equations
\begin{equation*}
\left\{
\begin{array}{l}
{\displaystyle 
\partial_t \rho + N^{b} \sum_{v\in\mc V} 
v \cdot \nabla F_0(\rho, \bs p)
\; =\; \Delta \rho \;,} \\
{\displaystyle
\partial_t p_j + N^{b} \sum_{v\in\mc V} v_j
\, v \cdot \nabla F_j(\rho, \bs p) \; =\; \Delta p_j\;,}
\end{array}
\right.
\end{equation*}
where $\rho$ stands for the density and $\bs p = (p_1, \dots, p_d)$
for the momentum. $F_0, \dots, F_d$ are thermodynamical quantities
determined by the ergodic properties of the dynamics.

Consider an initial profile given by $ (\rho, \bs p) = (\alpha, \bs
\beta) + N^{-b} (\varphi_0, \bs \varphi)$, where $(\alpha, \bs \beta)$
are appropriate constants. Expanding the solution of the previous
equations around $(\alpha, \bs \beta)$ and assuming that the first
component $\varphi_0$ does not depend on space, we obtain that the
momentum should evolve according to the incompressible Navier-Stokes
equation
\begin{equation*}
\left\{
\begin{array}{l}
{\displaystyle \vphantom{\Big\{}
\text{div }\bs \varphi =0 \;,} \\
{\displaystyle \vphantom{\Big\{}
\partial_t \varphi_\ell \;=\; A_0 \partial_\ell \varphi_\ell^2
\;+\; A_1 \varphi \cdot \nabla \varphi_\ell 
\;+\; A_2 \partial_\ell |\varphi |^2 \; 
+\; \Delta \varphi_\ell}\;,
\end{array}
\right.
\end{equation*}
for $1\le \ell\le d$, where $A_0$, $A_1$, $A_2$ are model-dependent
constants. This is the content of the main theorem of the article. We
prove that under an appropriate time scale the normalized empirical
measures associated to the momentum converge to the solution of the
above incompressible Navier-Stokes equation.

The proof relies on the relative entropy method introduced by Yau
\cite{y2}. We show that the entropy of the state of the process with
respect to a slowly varying parameter Gibbs state is small in a finite
time interval provided the solution of the incompressible
Navier-Stokes equation is smooth in this interval. 

To obtain such a bound on the entropy, we compute its time derivative
which can be expressed in terms of currents. A one block estimate,
which requires a polynomial bound on the spectral gap of the generator
of the process, permits to express the currents in terms of the
empirical density and momenta. The linear part of the functions of the
density and momenta cancel; while the second order terms can be
estimated by the entropy.  We obtain in this way a Gronwall inequality
for the relative entropy, which in turn give the required bound.

The article is organized as follows. In Section \ref{sec1} we
establish the notation and state the main results of the article.  In
Sections \ref{sec7} and \ref{sec2.1}, we examine the incompressible
limit of an asymmetric long range exclusion process. We state in this
simpler context some ergodic theorems needed in the proof of the
incompressible limit of the stochastic lattice gas. In Section
\ref{sec4} we prove the main result of the article, while in Section
\ref{gap} we prove a spectral gap, polynomial in the volume, for the
generator of a stochastic lattice gas restricted to a finite cube and
in Section \ref{equiv} we state an equivalence of ensembles for the
canonical measures of lattice gas models.

\section{Notation and Results}
\label{sec1}

Denote by $\bb T_N^d = \{0, \dots $, $N-1\}^d$ the $d$-dimensional
torus with $N^d$ points and let $ \mc V \subset \bb R^d $ be a finite
set of velocities $v=(v_1, \dots, v_d)$.  Assume that $\mc V$ is
invariant under reflexions and permutations of the coordinates:
\begin{equation*}
(v_1, \dots, v_{i-1}, -v_i, v_{i+1}, \dots, v_d) \quad \text{and}
\quad (v_{\sigma(1)}, \dots, v_{\sigma(d)}) 
\end{equation*}
belong to $\mc V$ for all $1\le i\le d$ and all permutations $\sigma$
of $\{1, \dots, d\}$ provided $(v_1, \dots, v_d)$ belongs to $\mc V$.

On each site of the discrete $d$-dimensional torus $\bb T_N^d$ at most
one particle for each velocity is allowed.  A configuration is denoted
by $\eta= \{\eta_x,~x\in \bb T_N^d\}$ where
$\eta_x=\{\eta(x,v),~v\in{\mc V}\}$ and $\eta(x,v)\in\{0,1\},~x\in\bb
T_N^d,~v\in{\mc V}$, is the number of particles with velocity $v$ at
$x$.  The set of particle configurations is $X_N=\left(\{0,1\}^{\mc
    V}\right)^{\bb T_N^d}$.

The dynamics consists of two parts: long range asymmetric random walks
with exclusion among particles of the same velocity and binary
collisions between particles of different velocities.  The first part
of the dynamics corresponds to the evolution of a mesoscopic
asymmetric simple exclusion process.  The jump law and the waiting
times are chosen so that the rate of jumping from site $x$ to site
$x+z$ for a particle with velocity $v$ is $p_N(z,v)$, where
\begin{equation*}
p_N(z,v) \;=\; \frac {A_M}{M^{d+2}} \big\{ 2  \; +\; \frac {1}{N^a} \,
q_M(z,v) \big\} \, \mb 1\{z\in\Lambda_M\}\;.
\end{equation*}
In this formula $a>0$ is a fixed parameter, $M$ is a function of $N$
to be chosen later, $\Lambda_M$ is the cube $\{-M, \dots, M\}^d$,
$A_M$ is given by
\begin{equation}
\label{g04}
\frac{A_M}{M^{d+2}} \sum_{z\in\Lambda_M} z_i z_j  \;=\; \delta_{i,j}
\end{equation}
for $1\le i,j\le d$, and $q_M(z,v)$ is any bounded non-negative
rate such that
\begin{equation*}
\frac {A_M}{M^{d+1}} \sum_{z\in\Lambda_M} q_M(z,v) z_i\;=\; v_i
\end{equation*}
for $1\le i\le d$ and $M\ge 1$. A possible choice is $q_M(z,v) =
M^{-1} (z\cdot v)$, where $u\cdot v$ stands for the inner product in
$\bb R^d$. Note that particles with velocity $v$ have mean
displacement $M^{-1} N^{-a} v$.

The generator ${\mc L}_N^{ex}$ of the random walk part of the dynamics
acts on local functions $f$ of the configuration space $X_N$ as
\begin{equation*}
(\mc L_N^{ex} f) (\eta) \;=\; 
\sum_{v \in \mc V} \sum_{\substack{x \in \bb T_N^d \\ z\in \Lambda_M}}
\eta(x,v) \, [1-\eta(x+z,v)] \, p_N(z,v) \,
[f(\eta^{x, x+z, v}) - f(\eta)]\;,
\end{equation*}
where
\begin{equation*}
\eta^{x,y,v}(z,w) \;=\;
\left\{
\begin{array}{ll}
\eta(y,v) & \text{if $w=v$ and $z=x$,} \\
\eta(x,v) & \text{if $w=v$ and $z=y$,} \\
\eta(z,w) & \text{otherwise.}
\end{array}
\right.
\end{equation*}

The collision part of the dynamics is described as follows. Denote by
${\mc Q}$ the set of all collisions which preserve momentum:
\begin{equation*}
\mc Q \;=\; \{(v,w,v',w') \in \mc V^4: v+w = v'+w'\}\;.
\end{equation*}
Particles of velocities $v$ and $w$ at the same site collide at rate
one and produce two particles of velocities $v'$ and $w'$ at that
site.  The generator $\mc L_N^c$ is therefore
\begin{equation*}
\mc L_N^c f(\eta) \;=\; \sum_{y \in \bb T_N^d}
\sum_{q \in \mc Q} p(y,q,\eta) \, [f(\eta^{y,q}) - f(\eta)]\;,
\end{equation*}
where the rate $p(y,q,\eta)$, $q=(v,w,v',w')$, is given by
\begin{equation*}
p(y,q,\eta) \;=\; \eta(y,v) \, \eta(y,w)\,
[1- \eta(y,v')]\, [1-\eta(y,w')]
\end{equation*}
and where the configuration $\eta^{y,q}$, $q= (v_0, v_1, v_2, v_3)$, 
after the collision is defined as
\begin{equation*}
\eta^{y,q}(z,u) \;=\; \left\{
\begin{array}{ll}
\eta(y,v_{j+2})  & \text{if $z=y$ and $u=v_j$ for some $0\le j\le 3$,}\\
\eta(z,u) & \text{otherwise,}
\end{array}
\right.
\end{equation*}
where the index of $v_{j+2}$ should be understood modulo $4$.

The generator $\mc L_N$ of the stochastic lattice gas we examine in
this article is the superposition of the exclusion dynamics with the
collisions just introduced:
\begin{equation*}
\mc L_N  \;=\; N^2\big\{ \mc L_N^{ex} \;+\; \mc L_N^c \big\} \;.
\end{equation*}
Note that time has been \emph{speeded up} diffusively.  Let $\{\eta
(t) : t\ge 0\}$ be the Markov process with generator $\mc L_N$ and
denote by $\{S_t^N : t\ge 0\}$ the semigroup associated to $\mc L_N$.

For a probability measure $\mu$ on $X_N$, denote by $\bb P_{\mu}$ the
measure on the path space $D(\bb R_+, X_N)$ induced by $\{\eta (t) :
t\ge 0\}$ and the initial measure $\mu$.  Expectation with respect to
$\bb P_{\mu}$ is denoted by $\bb E_{\mu}$.

\subsection{The invariant states}
\label{sec3.1}

For each configuration $\xi\in\{0,1\}^{\mc V}$, denote by $I_0 (\xi)$
the mass of $\xi$ and by $I_k (\xi)$, $k =1, \dots, d$, the momentum of
$\xi$: 
\begin{equation*}
I_0(\xi) \;=\; \sum_{v\in \mc V} \xi(v)\;, \quad
I_{k}(\xi) \;=\; \sum_{v\in \mc V} v_k \, \xi (v)\;.
\end{equation*}
Set $\bs I(\xi):=(I_0(\xi), \dots, I_d(\xi))$. Assume that the set of
velocities $\mc V$ is chosen in such a way that the unique quantities
conserved by the dynamics $\mc L_N$ are mass and momentum:
$\sum_{x\in\bb T^d_N}\bs I(\eta_x)$.

Two examples of sets of velocities with this property were proposed by
Esposito, Marra and Yau \cite{emy2}. In Model I, $\mc V= \{\pm
e_1,\dots, \pm e_d\}$, where $\{e_j, j=1,\dots ,d\}$ stands for the
canonical basis of $\bb R^d$. In Model II, $d=3$, $w$ is a root of
$w^4-6w^2 -1$ and $\mc V$ contains $(1,1, w)$, all reflections of this
vector and all permutations of the coordinates, performing a total of
24 vectors since $w\not = \pm 1$.

For each chemical potential $\bs \lambda = (\lambda_0, \lambda_1,
\dots, \lambda_d)$ in $\bb R^{d+1}$, denote by $m_{\bs \lambda}$ the
measure on $\{0,1\}^{\mc V}$ given by
\begin{equation*}
m_{\bs \lambda} (\xi) \;=\; \frac 1{Z(\bs \lambda)}
\exp\Big\{\bs \lambda\cdot \bs I(\xi)
\Big\}\;, 
\end{equation*}
where $Z(\bs \lambda)$ is a normalizing constant. Notice that $m_{\bs
  \lambda}$ is a product measure on $\{0,1\}^{\mc V}$, i.e., that the
variables $\{\xi(v)\,: v\in \mc V\}$ are independent under $m_{\bs
  \lambda}$.

Denote by $\mu^N_{\bs \lambda}$ the product measure on $(\{0,1\}^{\mc
  V})^{\bb T_N^d}$ with marginals given by
\begin{equation*}
\mu^N_{\bs \lambda} \{ \eta : \eta (x,\cdot) = \xi\}
\;=\; m_{\bs \lambda} (\xi)
\end{equation*}
for each $\xi$ in $\{0,1\}^{\mc V}$ and $x$ in $\bb T_N^d$. Notice
that $\{\eta (x,v)\, : x\in \bb T_N^d, v\in \mc V\}$ are independent
variables under $\mu^N_{\bs \lambda}$. 

For each $\bs \lambda$ in $\bb R^{d+1}$ a simple computation shows
that $\mu^N_{\bs\lambda}$ is an invariant state for the Markov process
with generator $\mc L_N$, that the generator $\mc L_N^c$ is symmetric
with respect to $\mu^N_{\bs \lambda}$ and that $\mc L^{ex}_N$ has an
adjoint $\mc L^{ex, *}_N$ in which $p_N(z,v)$ is replaced by
$p_N^*(z,v) = p_N(-z,v)$. In particular, if we denote by $\mc L^{ex,
  s}_N$, $\mc L^{ex, a}_N$ the symmetric and the anti-symmetric part
of $\mc L^{ex}_N$, we have that 
\begin{eqnarray*}
\!\!\!\!\!\!\!\!\!\!\!\!\!\!\! &&
(\mc L_N^{ex,s} f) (\eta) \;=\; \frac {2A_M}{M^{d+2}}
\sum_{v \in \mc V}  \sum_{\substack{x \in \bb T_N^d \\ z\in \Lambda_M}}
(T_{x,x+z,v} f)(\eta) \;,  \\
\!\!\!\!\!\!\!\!\!\!\!\!\!\!\! && \qquad
(\mc L_N^{ex,a} f) (\eta) \;=\; \frac {A_M}{M^{d+2} N^a}
\sum_{v \in \mc V}  \sum_{\substack{x \in \bb T_N^d \\ z\in
    \Lambda_M}} q_M(z,v) \, (T_{x,x+z,v} f)(\eta)  \;,
\end{eqnarray*}
where
\begin{equation*}
(T_{x,x+z,v} f)(\eta)\; =\; 
\eta(x,v) \, [1-\eta(x+z,v)] \, [f(\eta^{x, x+z, v}) - f(\eta)]\;.
\end{equation*}

The expectation under the invariant state $\mu^N_{\bs \lambda}$ of the
mass and momentum are given by
\begin{eqnarray*}
\rho(\bs \lambda) & := & E_{m_{\bs \lambda}}
[ I_0(\xi)] \;=\; \sum_{v\in {\mc V}}
\theta_v (\bs \lambda)  \\
p_k(\bs \lambda) & := & E_{m_{\bs \lambda}}
[ I_k (\xi) ] \;=\; \sum_{v\in {\mc V}} v_k \, 
\theta_v (\bs \lambda)\;.
\end{eqnarray*}
In this formula $\theta_v(\bs \lambda)$ denotes the expected value of
the density of particles with velocity $v$ under $m_{\bs \lambda}$:

\begin{equation}
\label{hf2}
\theta_v (\bs \lambda) \;:=\; E_{m_{\bs \lambda}} [\xi(v)]
\;=\; \frac{\exp\big\{ \lambda_0 + \sum_{k =1}^d
\lambda_k \, v_k \big\}}
{1+\exp\big\{ \lambda_0 + \sum_{k =1}^d
\lambda_k \, v_k \big\} }\;\cdot
\end{equation}

Denote by $(\rho,\bs p)(\bs \lambda):=(\rho(\bs \lambda) , p_1(\bs
\lambda), \dots, p_d(\bs \lambda))$ the map which associates the
chemical potential to the vector of density and momentum. Note that
$(\rho, \bs p)$ is the gradient of the strictly convex function $\log
Z(\bs \lambda)$. In particular, $(\rho, \bs p)$ is one to one. In
fact, it is possible to prove that $(\rho, \bs p)$ is a diffeomorphism
onto $\mf A\subset\bb R^{d+1}$, the interior of the convex envelope of
$\{\bs I(\xi), \xi\in\{0,1\}^{\mc V}\}$.  Denote by $\bs
\Lambda=(\Lambda_0,\dots,\Lambda_d) :\mf A \to \bb R^{d+1}$ the
inverse of $(\rho,\bs p)$. This correspondence permits to parameterize
the invariant states by the density and the momentum: for each $(\rho,
\bs p)$ in $\mf A$ we have a product measure $\nu^N_{\rho, \bs p} =
\mu^N_{\bs \Lambda(\rho, \bs p)}$ on $(\{0,1\}^{\mc V})^{\bb T_N^d}$.

\subsection{Spectral gap}
\label{sec3.3}

For $L\ge 1$ and a configuration $\eta$, let $\bs I^L(x) = (I^L_0(x),
\dots$, $I^L_d(x))$ be the average of the conserved quantities in a
cube of length $L$ centered at $x$:
\begin{equation}
\label{hf1}
\bs I^L(x) \;:=\; \bs I^L(x, \eta) \;=\; 
\frac{1}{|\Lambda_L|}\sum_{z \in x+ \Lambda_L} \bs I(\eta_z)\;.
\end{equation}
Let $\mf V_{L}$ be the set of all possible values of $\bs I^L(0)$ when
$\eta$ runs over $\big(\{0,1\}^{\mc V}\big)^{\Lambda_L}$.  Obviously
$\mf V_{L}$ is a finite subset of the convex envelope of $\{\bs
I(\xi): \xi\in\{0,1\}^{\mc V} \}$.  The set of configurations
$\big(\{0,1\}^{\mc V}\big)^{\Lambda_L}$ splits in invariant subsets:
For each $\bs i$ in $\mf V_{L}$, let
\begin{equation*}
\mc H_L (\bs i)\;:=\; \{\eta\in\big(\{0,1\}^{\mc V}\big)^{\Lambda_L} :
\bs I^{L}(0)=\bs i\}\; .
\end{equation*}
For each
$\bs i$ in $\mf V_L$, define the canonical measure $\nu_{\Lambda_L,\bs
  i}$ as the uniform probability measure on $\mc H_L(\bs i)$.

Denote by $\mc L_{\Lambda_M}$ the generator $\mc L_N$ restricted to
the cube $\Lambda_M$ without acceleration. More precisely, on the
state space $\big(\{0,1\}^{\mc V}\big)^{\Lambda_M}$ consider the
generator $\mc L_{\Lambda_M} = \mc L_{\Lambda_M}^{ex}+ \mc
L_{\Lambda_M}^c$, which acts on local functions $f:\big(\{0,1\}^{\mc
  V}\big)^{\Lambda_M}\mapsto \bb R$ as
\begin{eqnarray*}
(\mc L_{\Lambda_M}^{ex} f) (\eta) & = & 
\sum_{v \in \mc V} \sum_{\substack{x, y \in \Lambda_M \\|x-y|<M}}
\eta(x,v) \, [1-\eta(y,v)] \, p_N(y-x,v) \,
[f(\eta^{x, y, v}) - f(\eta)]\;,\\
(\mc L_{\Lambda_M}^c f)(\eta) & = & \sum_{y \in \Lambda_M}
\sum_{q \in \mc Q} p(y,q,\eta) \, [f(\eta^{y,q}) - f(\eta)]\;.
\end{eqnarray*}

Since the only conserved quantities are the total mass and momentum,
the process restricted to each component $\mc H_M(\bs i)$ is ergodic.
It has therefore a finite spectral gap: For each $\bs i$ in $\mf
V_{M}$, there exists a finite constant $C(M,\bs i)$ such that
\begin{equation*}
\< f; f\>_{\nu_{\Lambda_M,\bs i}} \;\le\; 
C(M,\bs i) \<f, (- \mc L_{\Lambda_M} f) \>_{\nu_{\Lambda_M,\bs i}}
\end{equation*}
for all functions $f$ in $L^2(\nu_{\Lambda_M,\bs i})$. Here and below
$\< f; f\>_{\nu}$ stands for the variance of $f$ with respect to a
measure $\nu$ and $\< \cdot , \cdot \>_{\nu}$ for the scalar product
in $L^2(\nu)$.

We shall assume that the inverse of the spectral gap increases
polynomially in the length of the cube: There exists $C_0>0$ and
$\kappa>0$ such that
\begin{equation}
\label{g12}
\max_{\bs i\in \mf V_M} C(M,\bs i) 
\;\le\; C_0 M^\kappa.
\end{equation}
We prove this hypothesis in Section \ref{gap} for Model I.

\subsection{Incompressible limit}
\label{sec3.2}
For $k=0, \dots, d$, denote by $\pi^{k,N}$ the empirical measure
associated to the $k$-th conserved quantity:
\begin{equation*}
\pi^{k,N} \;=\; N^{-d} \sum_{x\in\bb T_N^d} I_k(\eta_x) \delta_{x/N}\;,
\end{equation*} 
where $\delta_u$ stands for the Dirac measure concentrated on $u$.

Denote by $\<\pi^{k,N}, H\>$ the integral of a test function $H$ with
respect to an empirical measure $\pi^{k,N}$. To compute $\mc L_N
\<\pi^{k,N}, H\>$, note that $\mc L_N^c I_k(\eta_x)$ vanishes for
$k=0, \dots, d$ because the collision operators preserve local mass
and momentum. In particular, $\mc L_N \< \pi^{k,N}, H\> = N^2 \mc
L_N^{ex} \< \pi^{k,N},H\>$. To compute $\mc L^{ex}_N \<\pi^{k,N},
H\>$, consider separately the symmetric and the anti-symmetric part of
$\mc L^{ex}_N$. After two summations by parts and a Taylor expansion,
we obtain that
\begin{eqnarray*}
\!\!\!\!\!\!\!\!\!\!\!\!\! &&
N^2 \mc L_N^{ex,s} \<\pi^{0,N}, H\> \;=\; 
\< \pi^{0,N}, \Delta H\> \;+\; O(M/N) \;, \\
\!\!\!\!\!\!\!\!\!\!\!\!\! && \quad
N^2 \mc L_N^{ex,a} \<\pi^{0,N}, H\> \;=\;
\frac {N^{1-a}}{M} \frac 1{N^d} \sum_{j=1}^d \sum_{x\in\bb T_N^d}
(\partial_{u_j} H)(x/N) \tau_x W^M_{j} \; +\; O(N^{-a}) \; ,
\end{eqnarray*}
for every smooth function $H$. In this formula, $\Delta$ stands for
the Laplacian.  $\tau_x$ stands for the translation by $x$ on the
state space $ X_N$ so that $(\tau_x \eta )(y,v) = \eta (x+y, v)$ for
all $x$, $y$ in $\bb Z^d$, $v$ in $\mc V$, and $W^M_{j}$, $j=1, \dots
, d$, is the current given by
\begin{equation}
\label{g02}
W^M_{j}\;=\;  \frac {A_M}{M^{d+1}} \sum_{v\in\mc V}
\sum_{z\in\Lambda_M} q_M(z,v) \, z_j \,
\eta(0,v)\, \{ 1- \eta(z,v) \} \; .
\end{equation}

In the same way, for $1\le k\le d$, a long but simple computation
shows that
\begin{eqnarray*}
\!\!\!\!\!\!\!\!\!\!\!\!\! &&
N^2 \mc L_N^{ex, s} \< \pi^{k,N}, H\> \;=\; 
\< \pi^{k,N}, \Delta H\> \;+\; O(M/N) \;, \\
\!\!\!\!\!\!\!\!\!\!\!\!\! && \quad
N^2 \mc L_N^{ex, a} \<\pi^{k,N}, H\> \;=\;
\frac {N^{1-a}}{M} \frac 1{N^d} \sum_{j=1}^d \sum_{x\in\bb T_N^d}
(\partial_{u_j} H)(x/N) \tau_x W^M_{k,j} \; +\; O(N^{-a}) \; ,
\end{eqnarray*}
where $W^M_{k,j}$ is the current defined by
\begin{equation}
\label{g03}
W^M_{k,j} \;=\; \frac {A_M}{M^{d+1}} \sum_{v\in\mc V} v_k \, 
\sum_{z\in\Lambda_M} q_M(z,v) \, z_j \,
\eta(0,v)\, \{ 1- \eta(z,v) \} \; .
\end{equation}

The explicit formulas for $\mc L_N \< \pi^{k,N}, H\>$ permit to
predict the hydrodynamic behavior of the system under diffusive
scaling assuming local equilibrium.  By \eqref{g04}, the expectation
of the currents $W^M_{j}$, $W^M_{k,j}$ under the invariant state
$\mu^N_{\bs \lambda}$ are given by
\begin{equation*} 
E_{\mu^N_{\bs \lambda}}\big[ W^M_{j} \big] \;=\; 
\sum_{v\in\mc V} \chi (\theta_v(\bs \lambda)) \, v_j \ \;, \quad
E_{\mu^N_{\bs \lambda}}\big[ W^M_{k,j} \big] \;=\; 
\sum_{v\in\mc V} \chi (\theta_v(\bs \lambda)) \, v_j \, v_k \;.
\end{equation*}
In this formula and below, $\chi(a) = a(1-a)$.  In view of the
previous computation, if the conservation of local equilibrium holds,
the limiting equation in the diffusive regime is expected to be
\begin{equation}
\label{ff1}
\left\{
\begin{array}{l}
{\displaystyle 
\partial_t \rho + \frac {N^{1-a}}{M } \sum_{v\in\mc V} 
v \cdot \nabla \chi (\theta_v(\mb \Lambda(\rho, \bs p)))  
\; =\; \Delta \rho \;,} \\
{\displaystyle
\partial_t p_j + \frac {N^{1-a}}{M } \sum_{v\in\mc V} v_j
\, v \cdot \nabla \chi (\theta_v(\mb \Lambda(\rho, \bs p)))  
\; =\; \Delta p_j\;,}
\end{array}
\right.
\end{equation}
where $\nabla F$ stands for the gradient of $F$.

We turn now to the incompressible limit.  Note that $\theta_v(\bs 0) =
1/2$ for all $v$ and that
\begin{equation*}
\rho (\bs 0) = \frac{ |\mc V|}2 \;=: \; a_0\;, \quad p_k(\bs 0) \;=\;
\frac 12 \sum_{v\in\mc V} v_k  \;=\; 0\;,
\end{equation*}
where the last identity follows from the symmetry assumptions made on
$\mc V$. Therefore, $\bs \Lambda (a_0, \bs 0) = \bs 0$ and by
Taylor expansion, 
\begin{eqnarray*}
\!\!\!\!\!\!\!\!\!\!\!\!\!\! &&
\chi (\theta_v(\mb \Lambda(a_0 + \epsilon \varphi_0,
\epsilon \bs \varphi))) \;=\; \\
\!\!\!\!\!\!\!\!\!\!\!\!\!\! && \quad
\chi (1/2) \;-\;  \epsilon^2 \Big\{ \frac 14 \sum_{\ell =0}^d 
\partial_\ell \Lambda_0 (a_0,0) \varphi_\ell \;+\; \frac 14
\sum_{k=1}^d \sum_{\ell =0}^d v_k\,  
\partial_\ell \Lambda_k (a_0,0) \varphi_\ell \Big\}^2 \;+\;
O(\epsilon^3) 
\end{eqnarray*}
because $\chi'(1/2)=0$, $\partial_0 \theta_v (\bs 0) = (1/4)$,
$\partial_k \theta_v (\bs 0) = (1/4) v_k$, $1\le k\le d$.  Here
$\partial_\ell$ stands for the partial derivative with respect to the
$\ell$-th coordinate.  It follows from the previous explicit formulas
for $\partial_k \theta_v (\bs 0)$ that
\begin{equation*}
\partial_\ell \Lambda_0 (a_0, \bs 0) \;=\; 4 \delta_{0,\ell} 
|\mc V|^{-1}\;, \quad
\partial_\ell \Lambda_k (a_0, \bs 0) \;=\; 4 \delta_{k,\ell} 
\Big\{ \sum_{v\in\mc V} v_\ell^2 \Big\}^{-1} 
\end{equation*}
for $1\le k\le d$, $0\le \ell\le d$. In particular, $\chi
(\theta_v(\mb \Lambda(a_0 + \epsilon \varphi_0, \epsilon \bs
\varphi)))$ is equal to
\begin{equation*}
\chi (1/2) \;-\;  \epsilon^2 \Big\{  \frac{\varphi_0}{|\mc V|}  \;+\;
\sum_{k=1}^d \frac{ v_k\, \varphi_k}{\sum_{v\in\mc V} v_k^2}  \Big\}^2 
\;+\; O(\epsilon^3)
\end{equation*}
Due to the symmetry properties of $\mc V$, 
\begin{equation}
\label{g06}
\sum_{v\in\mc V} v_k \, v_j \; =\;  B \delta_{k,j} 
\end{equation}
where $B= \sum_{v\in\mc V} v_1^2$. The denominator in the expression
inside braces is thus equal to $B$.
 
To investigate the incompressible limit around $(a_0,\bs 0)$, fix
$b>0$ and assume that a solution of \eqref{ff1} has the form
$\rho (t,u) = a_0 + N^{-b} \varphi_0 (t,u)$, $p_k (t,u) = N^{-b}
\varphi_k (t,u)$. Then, to obtain a non-trivial limit we need to set
$M=N^{1-a-b}$ to obtain that $(\varphi_0,\bs \varphi)(t,u)$ is the
solution of
\begin{equation}
\label{g01}
\left\{
\begin{array}{l}
{\displaystyle 
\partial_t \varphi_0 \;=\; \sum_{v\in\mc V} 
v \cdot \nabla \Big\{  \frac{\varphi_0}{|\mc V|} \;+\; \frac 1B
\sum_{k=1}^d v_k \,\varphi_k \Big\}^2 \; 
+\; \Delta \varphi_0 \;,} \\
{\displaystyle
\partial_t \varphi_\ell \;=\; \sum_{v\in\mc V} 
v_\ell \, v \cdot \nabla \Big\{  \frac{\varphi_0}{|\mc V|} \;+\;
\frac 1B \sum_{k=1}^d v_k\,  \varphi_k \Big\}^2 \; 
+\; \Delta \varphi_\ell \;,}
\end{array}
\right.
\end{equation}
$1\le \ell \le d$.

To recover the Navier-Stokes equation, we need to introduce some
notation related to the velocity space $\mc V$.  Let $R_{k,\ell,m,n}
(v) \;=\; v_k \, v_\ell \, v_m \, v_n$. By the symmetry properties of
$\mc V$, if $k\not = \ell$, we have that
\begin{equation}
\label{g07}
\sum_{v\in\mc V} R_{k,\ell,m,n} (v) \;=\; \{ \delta_{m,k}
\delta_{n,\ell} + \delta_{m,\ell} \delta_{n,k}\} \, C\; ,
\end{equation}
where $C= \sum_{v\in\mc V} v_k^2 \, v_\ell^2 = \sum_{v\in\mc V} v_1^2
\, v_2^2$.  On the other hand,
\begin{equation}
\label{g08}
\sum_{v\in\mc V} R_{k,k,m,n} (v) \;=\;  \delta_{m,k}
\delta_{n,k} \, D \;+\; \delta_{m,n} \{ 1- \delta_{m,k}\} \, C\;,
\end{equation}
where $D = \sum_{v\in\mc V} v_k^4 = \sum_{v\in\mc V} v_1^4$.

Assume that $\varphi_0(0,u)$ is constant and that
$\text{div } \bs \varphi(0,\cdot) =0$. Since $\partial_t \varphi_0$ and
$\Delta \varphi_0$ vanish and since $\mc V$ is invariant by reflexion
around the origin, replacing $v$ by $-v$, the first equation in
\eqref{g01} can be rewritten as
\begin{equation*}
\sum_{k=1}^d \sum_{v\in\mc V} v_k \, 
v \cdot \nabla \varphi_k \;=\; 0\;.
\end{equation*}
By \eqref{g06}, this equation becomes
\begin{equation*}
\text{div }\bs \varphi =0\;.
\end{equation*}

The same argument permits to rewrite the second equations in
\eqref{g01} as
\begin{equation*}
\partial_t \varphi_\ell \;=\; B^{-2} \sum_{v\in\mc V} 
v_\ell \, v \cdot \nabla \Big\{
\sum_{k=1}^d v_k\,  \varphi_k \Big\}^2 \; 
+\; \Delta \varphi_\ell\;.
\end{equation*}
The first term on the right hand side of this expression is equal to
\begin{equation*}
2 B^{-2} \sum_{k,m,n =1}^d \varphi_k \, \partial_m \varphi_n 
\sum_{v\in\mc V} R_{k,\ell,m,n}(v)\;.
\end{equation*}
It follows from \eqref{g07}, \eqref{g08} and elementary algebra that
this expression is equal to $B^{-2}$ times
\begin{equation*}
(D-3C) \partial_\ell \varphi_\ell^2 \;+\; 2C \varphi \cdot \nabla
\varphi_\ell \;+\; C \partial_\ell |\varphi_\ell|^2
\end{equation*}
because $\text{div }\bs \varphi =0$. We recover in this way
Navier-Stokes equation
\begin{equation}
\label{g05}
\left\{
\begin{array}{l}
{\displaystyle \vphantom{\Big\{}
\text{div }\bs \varphi =0 \;,} \\
{\displaystyle \vphantom{\Big\{}
\partial_t \varphi_\ell \;=\; A_0 \partial_\ell \varphi_\ell^2
\;+\; A_1 \varphi \cdot \nabla \varphi_\ell 
\;+\; A_2 \partial_\ell |\varphi |^2 \; 
+\; \Delta \varphi_\ell}\;,
\end{array}
\right.
\end{equation}
where $A_0 = (D-3C)/B^2$, $A_1=2C/B^{-2}$ and $A_2=C/B^2$.  For Model
I we get $A_0=1$, $A_1=A_2=0$, while for Model II, $B=16 + 8w^2$, $C=8
+ 16 w^2$, $D=16 + 8w^4$ and $A_0$ vanishes because $w$ is chosen as a
root of $w^4-6w^2-1$.

\subsection{Statement of the result.}
\label{sec3.4}

Recall that $\kappa$ stands for the polynomial growth rate of the
spectral gap.  Assume that $b<a$,
\begin{equation}
\label{g09}
a\; +\; b \;>\; 1 - \frac 2{d+\kappa}\;,\quad
a\; +\; \Big(\frac {\kappa -2}{\kappa}\Big) b 
\;>\; 1 - \frac 2{\kappa}\;, \quad
a \;+\; \Big( 1 + \frac 2d \Big)b \;<\; 1\;.
\end{equation}
The first two displayed conditions are needed in the proof of the
one-block estimate, where the size of the cube cannot be too large.
The last condition appears in the replacement of expectations with
respect to canonical measures by expectations with respect to grand
canonical measures, where the volume $|\Lambda_M|$ has to be large.
It is easy to produce constants $a$, $b>0$ meeting the above
requirements. It is enough to choose first $0<a<1$, close enough to
$1$, and then to find $b$ small enough.


Let $\bs \varphi = (\varphi_1, \dots, \varphi_d): \bb T^d \to \bb R^d$
be a smooth divergence free vector field. Denote by $\bs \varphi(t)$
the solution of \eqref{g05} with initial condition $\bs \varphi$,
assumed to be smooth in a time interval $[0,T]$. Denote by $\nu^N_t$
the product measure on $X_N$ with chemical potential chosen so that
\begin{equation*}
E_{\nu^N_t} \big[ I_{0} (\eta_x) \big] \;=\; a_0 \; +\;
\frac{\varphi_0}{N^b}\; ,\quad
E_{\nu^N_t} \big[ I_{k} (\eta_x) \big] \;=\; \frac{\varphi_k(t,x/N)}{N^b}
\end{equation*}
for $1\le k\le d$ with $\varphi_0$ being a constant. This is possible
for $N$ large enough since $\bs\varphi$ is bounded and $\bs
\Lambda(a_0,\bs 0)=\bs 0$.  

For two probability measures $\mu$, $\nu$ on $X_N$, denote by $H_N(\mu
| \nu)$ the entropy of $\mu$ with respect to $\nu$:
\begin{equation*}
H_N(\mu | \nu) \;=\; \sup_f \Big\{ \int f d\mu \;-\; \log \int e^{f} d\nu
\Big\}\;,
\end{equation*}
where the supremum is carried over all bounded continuous functions on
$X_N$. We are now in a position to state the main theorem of this
article. 

\begin{theorem}
\label{gs01}
Assume conditions \eqref{g12} and \eqref{g09}.  Let $\bs \varphi =
(\varphi_1, \dots, \varphi_d): \bb T^d \to \bb R^d$ be a smooth
divergence free vector field. Denote by $\bs \varphi(t)$ the solution
of \eqref{g05} with initial condition $\bs \varphi$ and assume $\bs
\varphi(t,u)$ to be smooth in $[0,T] \times \bb T^d$ for some $T>0$.
Let $\{\mu^N: N\ge 1\}$ be a sequence of measures on $X_N$ such that
$H_N(\mu^N | \nu_0^N) = o(N^{d-2b})$.  Then, $H_N(\mu^N S_t^N |
\nu^N_t) = o(N^{d-2b})$ for $0\le t\le T$.
\end{theorem}

\begin{corollary}
\label{gs07}
Under the assumptions of Theorem \ref{gs01}, for every $0\le t\le T$
and every continuous function $F: \bb T^d\to\bb R$,
\begin{eqnarray*}
\!\!\!\!\!\!\!\!\!\!\!\!\!\! &&
\lim_{N\to\infty} \frac{N^b}{N^d} \sum_{x\in\bb T_N^d} F(x/N) 
\big\{I_0(\eta_x(t)) - a_0\big\} \;=\; \varphi_0 \int_{\bb T^d} F(u) 
\, du \;, \\
\!\!\!\!\!\!\!\!\!\!\!\!\!\! && \quad
\lim_{N\to\infty} \frac{N^b}{N^d} \sum_{x\in\bb T_N^d} F(x/N) 
I_k(\eta_x(t))  \;=\; \int_{\bb T^d} F(u) \varphi_k (t,u)
\, du
\end{eqnarray*}
in $L^1(\bb P_{\mu^N})$.
\end{corollary}

The corollary is an elementary consequence of the theorem and of the
entropy inequality.

\section{Mesoscopic asymmetric exclusion processes}
\label{sec7}

We start with a model with no velocities. The proof is simpler in this
context and the results stated will needed for the stochastic lattice
gas.  Denote by $\eta$ the configurations of the state space $\mc X_N
= \{0,1\}^{\bb T_N^d}$ so that $\eta (x)$ is either $0$ or $1$ if site
$x$ is vacant or not.  We consider a mesoscopic asymmetric exclusion
process on $\mc X_N$. This is the Markov process whose generator is
given by
\begin{equation*}
(L_N f) (\eta) \;=\; \sum_{\substack{x\in\bb T_N^d \\ z\in\bb Z^d}}
\eta(x) [1-\eta(x+z)] \, 
p_N(z) \, [f(\sigma^{x,x+z} \eta) - f(\eta)]\; ,
\end{equation*}
where,
\begin{equation*}
p_N(z) \;=\; \frac {1}{M^{d+2}} \big\{ 2 A_M \; +\; \frac 1 {N^a} \,
q(z) \big\} \, \mb 1\{z\in\Lambda_M\}\; .
\end{equation*}
In this formula $a>0$, $M$, $A_M$ are chosen as in the previous
section and $q(y) = \text{sign} (y \cdot v)$ for a fixed vector
$v\in\bb R^d$. On the other hand, $\sigma^{x,y}\eta$ is the
configuration obtained from $\eta$ by interchanging the occupation
variables $\eta(x)$, $\eta(y)$:
\begin{equation*}
(\sigma^{x,y} \eta) (z)\; =\;
\left\{
\begin{array}{ll}
\eta (z)  & \hbox{if $z\neq x,y$}\; , \\
\eta (y)  & \hbox{if $ z=x$}\; , \\
\eta (x)  & \hbox{if $ z=y$}\; .
\end{array}
\right.
\end{equation*}

For a probability measure $\mu$ on $\mc X_N$, $\bb P_{\mu}$ stands for
the measure on the path space $D(\bb R_+, \mc X_N)$ induced by the
Markov process with generator $L_N$ \emph {speeded up} by $N^2$ and
the initial measure $\mu$.  Expectation with respect to $\bb P_{\mu}$
is denoted by $\bb E_{\mu}$.  Denote by $\{S_t^N : t\ge 0\}$ the
semigroup associated to the generator $N^2 L_N$.

For $0\le \alpha\le 1$, denote by $\mu^N_\alpha$ the Bernoulli product
measure on $\mc X_N$ with density $\alpha$. An elementary computation
shows that $\mu^N_\alpha$ is an invariant state for the Markov process
with generator $L_N$. Moreover, the symmetric and the anti-symmetric
part of the generator $L_N$, respectively denoted by $L_N^s$, $L_N^a$,
are given by:
\begin{eqnarray*}
\!\!\!\!\!\!\!\!\!\!\!\!\!\!\! &&
(L_N^s f) (\eta) \;=\; \frac {2 A_M}{M^{d+2}}
\sum_{\substack{x\in\bb T_N^d \\ z\in\bb Z^d}}
\eta(x) [1-\eta(x+z)] \, [f(\sigma^{x,x+z} \eta) - f(\eta)] \;, \\
\!\!\!\!\!\!\!\!\!\!\!\!\!\!\! && \qquad
(L_N^a f) (\eta) \;=\; \frac {1}{M^{d+2} N^a}
\sum_{\substack{x\in\bb T_N^d \\ z\in\bb Z^d}} q(z)
\eta(x) [1-\eta(x+z)] \, [f(\sigma^{x,x+z} \eta) - f(\eta)] \;.
\end{eqnarray*}

We investigate in this and in the next section the incompressible
limit of this model.  Consider first the hydrodynamic behavior of the
process under diffusive scaling.  Denote by $\pi^N$ the empirical
measure associated to a configuration:
\begin{equation*}
\pi^N \;=\; \pi^N (\eta) \;=\; 
\frac 1{N^d} \sum_{x\in\bb T_N^d} \eta(x) \delta_{x/N}\;.
\end{equation*}

Denote by $\<\pi^{N}, H\>$ the integral of a test function $H$ with
respect to an empirical measure $\pi^{N}$. To compute $N^2 L_N
\<\pi^{N}, H\>$, we consider separately the symmetric and the
anti-symmetric part of the generator.  After two summations by parts
and a Taylor expansion, we obtain that
\begin{equation*}
N^2 L_N^s \<\pi^{N}, H\> \;=\; \< \pi^N, \Delta H\> \;+\; O(M/N) \;.
\end{equation*}
On the other hand, after a summation by parts, $N^2 L_N^a \<\pi^N,
H\>$ becomes
\begin{equation*}
\frac {N^{1-a}}{M} \frac 1{N^d} \sum_{x\in\bb T_N^d}
\sum_{j=1}^d (\partial_{u_j} H)(x/N) \tau_x W^M_j \; +\; O(N^{-a}) \; ,
\end{equation*}
where $\tau_x$ stands for the translation by $x$ on the state space
$\mc X_N$ so that $(\tau_x \eta )(y) = \eta (x+y)$ for all $x$, $y$ in
$\bb T_N^d$, and $\mb W^M = (W^M_1, \dots, W^M_d)$ is the current given by
\begin{equation}
\label{h05}
W^M_j\;=\;  \frac 1{M^{d+1}} \sum_{z\in\Lambda_M} q(z) \, z_j \,
\eta(0)\, \{ 1- \eta(z)\} \; .
\end{equation}

The expectation of the current under the invariant state
$\mu^N_\alpha$ is
\begin{equation}
\label{f03}
\alpha (1-\alpha) \frac 1{M^{d+1}}
\sum_{z\in\Lambda_M} q(z) \, z_j.
\end{equation}
Since $M= N^{1-a-b}$, the limiting equation in the diffusive regime is
therefore expected to be
\begin{equation*}
\partial_t \rho + N^b \, \mb \gamma \cdot \nabla
\rho(1-\rho) \; =\; \Delta \rho\;,
\end{equation*}
where $\gamma_j = \int_{[-1,1]^d} u_j \, q(u)\, du$.

To investigate the incompressible limit around density $1/2$, suppose
that a solution of the previous equation has the form $\rho (t,u) =
(1/2) + N^{-b} \varphi (t,u)$. An elementary computation shows that
\begin{equation*}
\partial_t \varphi \;=\; 
\mb \gamma \cdot \nabla \varphi^2 \;+\; \Delta \varphi\;.
\end{equation*}

Assume the following conditions on $a$ and $b$, which could
certainly be relaxed:
\begin{equation} 
\label{ass}
\frac d{d+2} < a+b\;, \quad a + \max 
\Big\{ 2, 1 + \frac 2d \Big\} b < 1\;,
\quad 2 \Big( 1 + \frac 2d \Big) b < 1\;.
\end{equation}
The first assumption, which forbids a large mesoscopic range $M$, is
used in the proof of the one-block estimate. The second and third
assumptions, which require a not too small range $M$, are used
throughout the proof to discard error terms.

By the same reasons of the previous section, there exist positive
constants $a$, $b$ satisfying these assumptions.

Fix a continuous function $\varphi_0 : \bb T^d \to \bb R$.  Denote by
$\varphi = \varphi(t,u)$ the solution of the nonlinear parabolic
equation
\begin{equation}
\label{f11}
\left\{
  \begin{array}{l}
\partial_t \varphi \;=\; \gamma \cdot \nabla \varphi^2
\;+\;  \Delta \varphi\;, \\
\varphi(0, \cdot) = \varphi_0(\cdot)\;.
  \end{array}
\right.
\end{equation}
For $t\ge 0$, let $\nu^N_t$ be the product measure on $\mc X_N$ with
marginals given by
\begin{equation*}
E_{\nu^N_t} \big[ \eta (x) \big] \;=\; \frac 12 \; +\;
\frac 1 {N^b} \varphi(t,x/N)\; .
\end{equation*}
This is possible for $N$ large enough because $\varphi$ is bounded.
Recall that we denote by $H_N(\nu | \mu )$ the relative
entropy of a probability measure $\nu$ with respect to $\mu$.

\begin{theorem}
\label{mt3}
Assume conditions \eqref{ass}. Fix a smooth function $\varphi_0: \bb
T^d \to \bb R$ and denote by $\varphi = \varphi(t,u)$ the solution of
\eqref{f11} with initial condition $\varphi_0$. Assume $\varphi$ to be
smooth in the layer $[0,T] \times \bb T^d$. Let $\{\mu^N: N\ge 1\}$ be
a sequence of measures on $\mc X_N$ such that $H_N(\mu^N | \nu^N_{0})
= o(N^{d-2b})$.  Then, $H_N(\mu^N S_t^N | \nu^N_t) = o(N^{d-2b})$ for
all $0\le t\le T$.
\end{theorem}

Fix two bounded functions $\varphi_i: \bb T^d \to \bb R$, $i=1, 2$,
and denote by $\nu^{N,i}$ the product measures associated to the
density profile $(1/2) + N^{-b} \varphi_i$. A second order Taylor
expansion shows that
\begin{equation*}
H_N(\nu^{N,2} | \nu^{N,1}) \; = \; O( N^{d-2b} )\;.
\end{equation*}
The assumption on the entropy formulated in the theorem permits
therefore to distinguish between $N^{-b}$-perturbations of a constant
density profile.

A law of large numbers for the corrected empirical measure follows
from this result.  For a configuration $\eta$, denote by $\Pi^N
(\eta)$ the corrected empirical measure defined by
\begin{equation*}
\Pi^N \;=\; \Pi^N (\eta) \;=\;  \frac {N^b}{N^{d}} \sum_{x\in\bb T_N^d}
\{ \eta (x) - 1/2\} \, \delta_{x/N}
\end{equation*}
considered as an element of $\mc M(\bb T^d)$, the space of Radon
measures on $\bb T^d$ endowed with the weak topology. For $t\ge 0$, let
$\Pi^N_t = \Pi^N (\eta_t)$. 

\begin{corollary}
\label{mt4}
Under the assumptions of Theorem \ref{mt3}, for every $0\le t\le T$
and every continuous function $F: \bb T^d\to\bb R$,
\begin{eqnarray*}
\lim_{N\to\infty} \< \Pi^N_t,  F \>
\;=\;  \int_{\bb T^d} \varphi (t,u) F(u) \, du 
\end{eqnarray*}
in $L^1(\bb P_{\mu^N})$.
\end{corollary}

The corollary is an elementary consequence of Theorem \ref{mt3} and
the entropy inequality.

\section{Incompressible limit of mesoscopic exclusion processes.}
\label{sec2.1}

We prove in this section Theorem \ref{mt3}. Fix a smooth function
$\varphi_0: \bb T^d \to \bb R$ and denote by $\varphi = \varphi(t,u)$
the solution of \eqref{f11} with initial condition $\varphi_0$,
supposed to be smooth in the time interval $[0,T]$. Let $\{\mu^N: N\ge
1\}$ be a sequence of measures on $\mc X_N$ satisfying the assumptions
of Theorem \ref{mt3}.

\subsection{Entropy, Dirichlet form and Ergodic constants}
\label{secb1}

An elementary computation shows that the entropy $H_N(\mu^N |
\mu^N_{1/2})$ of $\mu^N$ with respect to $\mu^N_{1/2}$ is of order
$N^{d-2b}$. Indeed, by the explicit formula for the entropy and by the
entropy inequality,
\begin{equation*}
H_N(\mu^N | \mu^N_{1/2}) \;\le\; \big( 1 + \frac 1A \big)
H_N(\mu^N | \nu_0^N) \; +\; \frac 1A \log \int \big(
\frac{d \nu_0^N}{d \mu^N_{1/2}}\big)^{1+A} d\mu^N_{1/2}
\end{equation*}
for all $A >0$. A Taylor expansion shows that the second term on
the right hand side is of order $N^{d-2b}$. In particular
\begin{equation}
\label{f02}
N^{2b-d} H_N(\mu^N | \mu^N_{1/2}) \;\le\; C_0
\end{equation}
for some finite constant $C_0$ depending only on $\varphi_0$.

Let $f_t^N$ be the Radon-Nikodym derivative $d \mu^N S_t^N / d
\mu^N_{1/2}$ so that
\begin{equation*}
  \partial_t f_t^N \;=\; N^2 L_N^* f_t^N\;,
\end{equation*}
where $L_N^*$ stands for the adjoint of $L_N$ in $L^2(\mu^N_{1/2})$.
It follows from \eqref{f02} and a well known estimate on the entropy
production (cf. \cite{kl}, Section V.2) that 
\begin{equation}
\label{f07}
\frac{N^{2b}}{N^d} H_N(\mu^N S_t^N | \mu^N_{1/2}) 
\;+\; \frac{N^{2b}}{N^d}
\int_0^t D_N(\mu^N_{1/2}, f_s^N) ds \;\le\; C_0
\end{equation}
for all $N\ge0$ and $t\ge 0$. In this formula, $D_N$ stands for
the Dirichlet form defined as
\begin{equation*}
D_N(\mu^N_{1/2}, f) \;=\; N^2 \< - L_N \sqrt{f}, 
\sqrt{f}\>_{\mu^N_{1/2}}\; ,
\end{equation*}
where $\< \cdot, \cdot \>_{\mu^N_{1/2}}$ is the scalar product in
$L^2(\mu^N_{1/2})$. An elementary computation shows that
\begin{equation*}
D_N (\mu^N_{1/2}, f) \;=\; \frac {A_M N^2}{M^{d+2}} \sum_{|x-y|\le M}
\< \{ \nabla^{x,y} \sqrt{f} \}^2 \>_{\mu^N_{1/2}}\; ,
\end{equation*}
where
\begin{equation*}
(\nabla^{x,y} g)(\eta) \;=\; g(\sigma ^{x,y} \eta) - g(\eta)\;,
\end{equation*}
and that $D_N$ is a convex, lower semicontinuous functional.

Let $L_{\Lambda_M}$ be the symmetric part of the generator $L_N$
restricted to the cube $\Lambda_M$:
\begin{equation}
\label{f04} 
(L_{\Lambda_M} f) (\eta) \;=\; \frac {2A_M}{M^{d+2}}
\sum_{\substack{x,y\in\Lambda_M \\ |x-y|\le M}} \eta(x) [1-\eta(y)]
\, [f(\sigma^{x,y} \eta) - f(\eta)]\; ,
\end{equation}
and denote by $\mu_{\Lambda_M, K}$, $0\le K\le |\Lambda_M|$, the
canonical measure on $\{0, 1\}^{\Lambda_M}$ concentrated on the
hyperplane with $K$ particles. In the case of the exclusion process,
$\mu_{\Lambda_M, K}$ is just the uniform measure over all
configurations of $\{0, 1\}^{\Lambda_M}$ with $K$ particles. Denote by
$D_{\Lambda_M}$ the Dirichlet form associated to $L_{ \Lambda_M}$:
\begin{equation*}
D_{\Lambda_M} (\mu, f) \;=\; \frac {A_M}{M^{d+2}}
\sum_{\substack{x,y \in \Lambda_M \\ |x-y|\le M}}
\< \{ \nabla^{x,y} \sqrt{f} \}^2 \>_{\mu}\; ,
\end{equation*}
where $\mu$ stands either for the marginal on $\Lambda_M$ of the grand
canonical measure $\mu^N_{1/2}$ or for a canonical measure
$\mu_{\Lambda_M, K}$.

By comparing the Dirichlet form $D_{\Lambda_M}$ with the
Bernoulli-Laplace Dirichlet form, in which all jumps are allowed with
rate $|\Lambda_M|^{-1}$ and which is known to have a spectral gap of
order $1$ (cf. \cite{q}), we can prove that the spectral gap of $D_{
  \Lambda_M}$ is of order $M^{-2}$.

\subsection{The relative entropy method}

The proof of Theorem \ref{mt3} is based on the relative entropy method
introduced by Yau \cite{y2}. Let $\psi_t^N = d \nu^N_t/ d\mu^N_{1/2}$.
It follows from the explicit formulas for the product measure
$\nu^N_{t}$ that
\begin{equation*}
\log \psi_t^N = \sum_{x\in \bb T_N^d} \log 
\frac { [(1/2) + N^{-b} \varphi(t,x/N)] } 
{[(1/2) - N^{-b} \varphi(t,x/N)]} \, 
\eta (x) \; + \; \sum_{x\in \bb T_N^d} 
\log \big\{ 1 - 2 N^{-b} \varphi(t,x/N) \big\}  \;.
\end{equation*} 

Let $H_N(t) = N^{2b - d} H_N(\mu^N S_t^N | \nu^N_t)$ and recall that
we denote the Radon-Niko\-dym derivative $d\mu^N S_t^N / d\mu^N_{1/2}$
by $f^N_t$. With the notation just introduced, we have that
\begin{equation*}
H_N(t) \;=\; N^{2b - d} \int f_t^N \log \frac{f_t^N}{\psi_t^N} 
\, d\mu^N_{1/2}\;.
\end{equation*}
Theorem \ref{mt3} follows from Gronwall lemma and the following
estimate.

\begin{proposition}
\label{h01}
Fix a sequence of measures $\{\mu^N : N\ge 1\}$ satisfying the
assumptions of Theorem \ref{mt3}.  There exists $\gamma >0$ such that
\begin{equation*}
H_N(t) \;\le\; \gamma \int_0^t H_N(s) \, ds \;+\; o_N(1)
\end{equation*}
for all $t\le T$.
\end{proposition}

The proof of Proposition \ref{h01} is divided in several steps. We
begin with a well known upper bound for the entropy production (see
e.g. \cite{kl}, Lemma 6.1.4).
\begin{equation}
\label{hf01}
\frac{d}{dt} H_N(t) \;\le\; N^{2b-d}
\int f_t^N \frac{( N^2 L_N^* - \partial_t) \psi_t^N}{\psi_t^N} 
\, d\mu^N_{1/2}\;.
\end{equation}

A long and tedious computation gives that
$(\psi_t^N)^{-1} (N^2 L_N^* - \partial_t) \psi_t^N$ is equal to
\begin{eqnarray}
\label{hf02}
\!\!\!\!\!\!\!\!\!\!\!\!\!\!\! && \nonumber
\frac{4}{N^b} \sum_{x \in \bb T_N^d} (\Delta
\varphi)(t,x/N) \{ \eta (x) - 1/2\}  \\
\!\!\!\!\!\!\!\!\!\!\!\!\!\!\! && \quad
+\; \frac{16}{N^{2b}} \sum_{i,j=1}^d  
\sum_{x \in \bb T_N^d} (\partial_{u_i} \varphi) (t,x/N)
(\partial_{u_j} \varphi ) (t,x/N) \tau_x V_{i,j}^{M} (\eta) \\
\!\!\!\!\!\!\!\!\!\!\!\!\!\!\! && \nonumber \qquad
\;+\; 4 ( 1 + \varepsilon_N) \sum_{j = 1}^d
\sum_{x \in \bb T_N^d} (\partial_{u_j} \varphi) (t,x/N) 
\tau_x W_{j}^{*,M} (\eta) \\
\!\!\!\!\!\!\!\!\!\!\!\!\!\!\! && \nonumber \qquad\quad
-\; \frac 4{N^{b}} \sum_{x \in \bb T_N^d} (\partial_t \varphi) (t,x/N)
\big\{ \eta (x) - (1/2) - N^{-b} \varphi (t,x/N) \big\}\; 
+ \; o(N^{d-2b})\; .
\end{eqnarray}
In this formula, $W_{j}^{*,M}$ and $V_{i,j}^{M}$ stand for
\begin{eqnarray*}
\!\!\!\!\!\!\!\!\!\!\!\!\!\!\! && 
W_{j}^{*,M} (\eta) \;=\; \frac {1}{M^{d+1}} 
\sum_{z\in\Lambda_M} q(-z) \, z_j \, \eta (0) \, 
\{ 1- \eta (z)\} \;, \\
\!\!\!\!\!\!\!\!\!\!\!\!\!\!\! && 
\quad 
V_{i,j}^{M} (\eta) \;=\; \frac {A_M}{M^{d+2}} 
\sum_{z\in \Lambda_M} z_i \, z_j \, \eta (0) \, 
\{1-\eta (z)\}\;. 
\end{eqnarray*}
We used the inequalities $a>b$, $a+ 2b <1$, which follow from
assumptions \eqref{ass}, to estimate several terms in the above
computation by $o(N^{d-2b})$. The expression $(1 + \varepsilon_N)
(\partial_{u_j} \varphi)$ in the third line stands for
$(\partial_{u_j} \varphi) \{ 1 - 4 N^{-2b} \varphi^2\}^{-1}$. Keep in
mind that $\varepsilon_N$ is of order $N^{-2b}$.

If we replace $\eta (x) - 1/2$ in the first term of \eqref{hf02} by
$\eta (x) - (1/2) - N^{-b} \varphi (t,x/N)$ and $V_{i,j}^{M}$ by
$V_{i,j}^{M} - (1/4) \delta_{i,j}$, as $N\uparrow\infty$, the
expressions added multiplied by $(1/4) N^{2b-d}$ converge to
\begin{equation*}
\int_{\bb T^d} (\Delta \varphi ) \, \varphi \;+\;
\sum_{j=1}^d \int_{\bb T^d} (\partial_{u_j} \varphi )^2 \;=\; 0\;.
\end{equation*}

Therefore, in view of \eqref{hf01}, \eqref{hf02}, the time derivative
of the renormalized entropy $H_N(t)$ is bounded above by
\begin{eqnarray}
\label{hf03}
\!\!\!\!\!\!\!\!\!\!\!\!\!\!\! && 
\bb E_{\mu^N} \Big[
\frac {4 (1 + \varepsilon_N) N^{2b}} {N^d}  \sum_{j=1}^d 
\sum_{x \in \bb T_N^d}
(\partial_{u_j} \varphi ) (t,x/N) \tau_x W_{j}^{*,M} (\eta_t) \Big] \\
\!\!\!\!\!\!\!\!\!\!\!\!\!\!\! && \nonumber
+\; \bb E_{\mu^N} \Big[ \frac{16}{N^{d}} \sum_{i,j=1}^d   
\sum_{x \in \bb T_N^d} (\partial_{u_i} \varphi) (t,x/N)
(\partial_{u_j} \varphi ) (t,x/N) 
\tau_x  \hat V_{i,j}^{M} (\eta_t) \Big] \\
\!\!\!\!\!\!\!\!\!\!\!\!\!\!\! && \nonumber
+\; \bb E_{\mu^N} \Big[ 
\frac {4 N^b}{N^{d}} 
\sum_{x\in\bb T^d}( \Delta \varphi - \partial_t
\varphi) (t,x/N)
\Big\{ \eta_t(x) - (1/2) - N^{-b} \varphi (t,x/N) \Big\}\Big] \; 
+ \; o_N(1)\; ,
\end{eqnarray}
where $\hat V_{i,j}^{M}(\eta) = V_{i,j}^{M}(\eta) - (1/4)
\delta_{i,j}$.

We now use the ergodicity to replace the functions $W_{j}^{*,M}$ and
$\hat V_{i,j}^{M}$ by their projections on the conserved quantity over
mesoscopic cubes.

\subsection{One block estimate}
\label{ss1}

Recall from the end of Subsection \ref{secb1} that $\mu_{\Lambda_M,
  K}$ stands for the uniform measure over all configurations of $\{0,
1\}^{\Lambda_M}$ with $K$ particles.  For $1\le j\le d$, denote by
$F_j(K/|\Lambda_M|)$ the expected value of the current $W_j^{*,M}$
with respect to $\mu_{\Lambda_M,K}$. An elementary computation shows
that
$$
F_{j} (\beta) \;=\; E_{\mu_{\Lambda_M,K}} [W_j^{*,M}]
\;=\; - \gamma_j^M \, \Big \{ 1 +
\frac {1} {|\Lambda_M| -1} \Big\}\, \beta (1-\beta) \;,
$$
provided $\beta = K/|\Lambda_M|$ and $\gamma_j^M = M^{-(d+1)}
\sum_{z\in\Lambda_M} z_j q(z) = \gamma_j + O(M^{-1})$.

For a positive integer $\ell\ge 1$, let $\eta^\ell (x)$ be the average
number of particles in a cube of size $\ell$ around $x$:
$$
\eta^\ell (x) \;=\; \frac 1{|\Lambda_\ell |} \sum_{y\in x + \Lambda_\ell}
\eta(y)\; .
$$

For $M\ge 1$, $1\le j \le d$, let
\begin{equation*}
V_{j,M} \;=\; W_j^{*,M} \;-\; F_{j} (\eta^M (0)) \;.
\end{equation*}

\begin{lemma}
\label{s1}
For every $t>0$, $1\le j\le d$ and continuous function $G:\bb
T^d\to\bb R$,
$$
\lim_{N\to\infty} \bb E_{\mu^N} \Big[ \,\Big\vert
\int_0^t ds\, \frac {N^{2b}}{N^d} \sum_{x\in\bb T_N^d} G(x/N)
\tau_x V_{j,M} (\eta_s) \Big\vert \, \Big] \;=\; 0\;.
$$
\end{lemma}

\begin{proof}
By the entropy inequality and Jensen inequality, the expectation
appearing in the statement of the lemma is bounded above by
\begin{equation*}
\frac{N^{2b}}{A N^d}H_N(\mu^N | \mu^N_{1/2}) \;+\; \frac {N^{2b}}{AN^d}
\log \bb E_{\mu^N_{1/2}} \Big[ \exp A \Big\vert
\int_0^t ds\,  \sum_{x\in\bb T_N^d} G(x/N)
\tau_x V_{j,M} (\eta_s) \Big\vert \, \Big]
\end{equation*}
for every $A>0$. In view of \eqref{f02}, to prove the lemma it is
enough to show that the second term vanishes, as $N\uparrow\infty$,
for any $A>0$. Since $e^{|x|} \le e^x + e^{-x}$, it is enough to
estimate the previous expectation without the absolute value.

By Feynman-Kac formula and by the variational formula for the largest
eigenvalue of an operator, the second term without the absolute value
is bounded above by
\begin{equation}
\label{f05}
\frac{tN^{2b}}{AN^d} \sup_{f} \Big\{ \sum_{x\in\bb T_N^d} AG(x/N)
\int  \tau_x V_{j,M} \, f \,
d\mu^N_{1/2} \;-\; D_N(\mu^N_{1/2}, f) \Big\}\;,
\end{equation}
where the supremum is carried over all density functions $f$ with
respect to $\mu^N_{1/2}$.

Since the measure
$\mu^N_{1/2}$ is translation invariant and since $V_{j,M}$ depends on
the configuration only through $\{\eta (z) \, : z \in \Lambda_M\}$,
\begin{equation*}
\int  (\tau_x V_{j,M}) \, f \, d\mu^N_{1/2} \;=\;
\int  V_{j,M} \, (\tau_{-x} f) \, d \mu^N_{1/2} \;=\;
\int  V_{j,M}  \, f_{x,M} \, d \mu^N_{1/2}\; ,
\end{equation*}
where $f_{x,M} = E_{\mu^N_{1/2}} [ \tau_{-x} f | \eta(z) \, , z\in
\Lambda_M]$.

On the other hand, by convexity of the Dirichlet form and by
translation invariance of $\mu^N_{1/2}$, for any $x$, $y$ in
$\Lambda_M$ such that $|x-y|\le M$,
\begin{equation*}
\< \{ \nabla^{x,y} \sqrt{f_{z,M}} \}^2 \>_{\mu^N_{1/2}} \;\le\;
\< \{ \nabla^{x,y} \sqrt{\tau_{-z} f } \}^2 \>_{\mu^N_{1/2}}
\;=\; \< \{ \nabla^{x+z,y+z} \sqrt{f } \}^2 \>_{\mu^N_{1/2}}\; .
\end{equation*}
Therefore, summing over $x$, $y$ in $\Lambda_M$, $|x-y|\le M$ and in
$z$ in $\bb T_N^d$, we obtain that
\begin{equation*}
\sum_{z\in\bb T_N^d} \sum_{\substack{x,y \in \Lambda_M \\ |x-y|\le M}}
\< \{ \nabla^{x,y} \sqrt{f_{z,M}} \}^2 \>_{\mu^N_{1/2}} \;\le\;
C_0 M^d \sum_{|x-y|\le M} \< \{ \nabla^{x,y} \sqrt{f} \}^2 \>_{\mu^N_{1/2}}
\end{equation*}
for some universal constant $C_0$.

Recall the definition of the Dirichlet forms $D_N$ and
$D_{M,\Lambda_M}$ introduced above.  It follows from the previous
estimates that the expression \eqref{f05} is bounded
above by
\begin{equation}
\label{f09}
\frac{tN^{2b}}{AN^d}\sum_{x\in\bb T_N^d} 
\sup_f \Big\{ AG(x/N) \int  V_{j,M}\,  f
d\mu^N_{1/2} \;-\; \frac{C_0 N^2}{M^d} D_{M, \Lambda_M}
(\mu^N_{1/2}, f) \Big\}
\end{equation}
where the supremum is carried over all densities $f$ with respect to
the marginal of $\mu^N_{1/2}$ on the cube $\Lambda_M$.

In particular, by projecting the density over each hyperplane with a
fixed total number of particles and recalling the perturbation theorem
on the largest eigenvalue of a symmetric operator (Theorem 1.1 of
Appendix 3 in \cite{kl}), in view of \eqref{f09}, we obtain that
\eqref{f05} is less than or equal to
\begin{equation*}
C(\gamma) A t \Vert G \Vert_\infty^2 \frac{N^{2b} M^d}{N^2}
\< (- L_{M,\Lambda_M})^{-1} V_{j,M} , V_{j,M} \>_{\mu^N_{1/2}}
\end{equation*}
for some finite constant $C (\gamma)$ depending only on $\gamma$.
Here we need the assumption that $M^{d+2} \ll N^2$ to be allowed to
apply the Rayleigh expansion.

Since the generator $L_{M,\Lambda_M}$ has a spectral gap of order
$M^{-2}$, $\< (- L_{M,\Lambda_M})^{-1} V_{j,M} $, $V_{j,M}
\>_{\mu^N_{1/2}}$ is bounded by $C_0 M^2 \< V_{j,M} ; V_{j,M}
\>_{\mu^N_{1/2}}$, which is less than or equal to $C_0$ $M^{2-d}$.
Thus, \eqref{f05} is bounded by $C_0 A t \Vert G \Vert_\infty^2
N^{-2a}$ because $M=N^{1-a-b}$. This concludes the proof of the lemma.
\end{proof}

For $1\le i, j\le d$, let
\begin{equation*}
F_{i,j} (\beta) \;=\; E_{\mu_{\Lambda_M,K}} [\hat V_{i,j}^{M}]
\;=\;  \beta (1-\beta) \, \Big \{ 1 +
\frac {1} {|\Lambda_M| -1} \Big\} \delta_{i,j} - (1/4) \delta_{i,j}\;,
\end{equation*}
with the same convention that $\beta = K/|\Lambda_M|$. Let $w_{i,j}^{M}
(\eta) = \hat V_{i,j}^{M} - F_{i,j} (\eta^M(0))$.  The arguments of
the proof of Lemma \ref{s1} shows that for every $t>0$, $1\le i,j\le
d$ and continuous function $G:\bb T^d\to\bb R$,
\begin{equation}
\label{hf05}
\lim_{N\to\infty} \bb E_{\mu^N} \Big[ \,\Big\vert
\int_0^t ds\, \frac {1}{N^d} \sum_{x\in\bb T_N^d} G(x/N)
\tau_x w_{i,j}^M (\eta_s) \Big\vert \, \Big] \;=\; 0\;.
\end{equation}
The arguments are even simpler due to the absence of the factor
$N^{2b}$ multiplying the sum.

By Lemma \ref{s1} and \eqref{hf05}, integrating in time \eqref{hf03},
we obtain that the entropy $H_N(t)$ is less than or equal to
\begin{eqnarray}
\label{hf04}
\!\!\!\!\!\!\!\!\!\!\!\!\!\!\! && \nonumber
\frac {4 (1 + \varepsilon_N) N^{2b}} {N^d}  
\int_0^t ds \, \bb E_{\mu^N} \Big[ \sum_{j=1}^d 
\sum_{x \in \bb T_N^d} (\partial_{u_j} \varphi ) (s,x/N) 
F_j (\eta^M_s(x)) \Big] \\
\!\!\!\!\!\!\!\!\!\!\!\!\!\!\! && 
+\; \frac{16}{N^{d}}  \int_0^t ds \,  \bb E_{\mu^N} \Big[ \sum_{i,j=1}^d   
\sum_{x \in \bb T_N^d} (\partial_{u_i} \varphi) (s,x/N)
(\partial_{u_j} \varphi ) (s,x/N) F_{i,j} (\eta^M_s(x))\Big] \\
\!\!\!\!\!\!\!\!\!\!\!\!\!\!\! && \nonumber
+\; \frac {4 N^b}{N^{d}}  \int_0^t ds \,  \bb E_{\mu^N} \Big[ 
\sum_{x\in\bb T^d}( \Delta \varphi - \partial_s \varphi) 
(s,x/N)
\Big\{ \eta_s(x) - (1/2) - N^{-b} \varphi (s,x/N) \Big\}\Big] 
\end{eqnarray}
plus an error term of order $o_N(1)$ for every $t\le T$.

Recall that $\chi (a) = a(1-a)$.  Since $N^{2b} \ll M^d$, we may
replace in the previous formula $F_j (\eta^M_s(x))$ by $-\gamma_j^M
\chi (\eta^M_s(x))$ and $F_{i,j} (\eta^M_s(x))$ by $\{ \chi
(\eta^M_s(x)) - \chi (1/2)\} \delta_{i,j}$. Moreover, since $N^{-d}
\sum_{x \in \bb T_N^d} (\partial_{u_j} \varphi ) (t,x/N)$ is of order
$N^{-1}$ and since $b<1/2$, we may further replace $\chi (\eta_t^M
(x))$ by $\chi (\eta_t^M (x)) - \chi (1/2)$ in the first term.
Finally, since for a smooth function $G$, $M^{-d} \sum_{y: |y-x|\le M}
[G(y/N) - G(x/N)]$ is of order $(M/N)^2$ and since $M^2 N^{b-2}$
vanishes as $N\uparrow\infty$, we may replace $\eta_s(x)$ by
$\eta^M_s(x)$ in the third term. After all these replacements and
since $\chi(b) - \chi(1/2) = -[b-(1/2)]^2$, \eqref{hf04} becomes
\begin{eqnarray}
\label{hf06}
\!\!\!\!\!\!\!\!\!\!\!\!\!\!\! && \nonumber
\frac {4 ( 1 + \varepsilon_N) N^{2b}} {N^d}
\int_0^t ds \, \bb E_{\mu^N} \Big[ \sum_{j=1}^d 
\sum_{x \in \bb T_N^d} \gamma_j^M (\partial_{u_j} \varphi ) (s,x/N) 
\big\{ \eta_s^M (x) - 1/2 \big\}^2 \Big] \\
\!\!\!\!\!\!\!\!\!\!\!\!\!\!\! && 
- \; \frac{16}{N^{d}} \int_0^t ds \,  \bb E_{\mu^N} \Big[ \sum_{j=1}^d   
\sum_{x \in \bb T_N^d} (\partial_{u_i} \varphi)^2 (s,x/N)
\big\{  \eta_s^M (x) - 1/2 \big\}^2\Big] \\
\!\!\!\!\!\!\!\!\!\!\!\!\!\!\! && \nonumber
+\; \frac {4 N^b}{N^{d}}  \int_0^t ds \,  \bb E_{\mu^N} \Big[ 
\sum_{x\in\bb T^d}( \Delta \varphi - \partial_s \varphi) (s,x/N)
\Big\{ \eta^M_s(x) - (1/2) - N^{-b} \varphi (s,x/N) \Big\}\Big] \;.
\end{eqnarray}

The second line of the previous formula is easy to estimate.  One can
argue that it is negative or one can add $N^{-b} \varphi (s,x/N)$
inside the braces and apply Lemma \ref{h09} below.  The first term in
\eqref{hf06} without the factor $(1 + \varepsilon_N)$ can be written
as
\begin{eqnarray*}
\!\!\!\!\!\!\!\!\!\!\!\!\!\!\! && 
\frac {4 N^{2b}} {N^d}
\int_0^t ds \, \bb E_{\mu^N} \Big[ \sum_{j=1}^d 
\sum_{x \in \bb T_N^d} \gamma_j^M (\partial_{u_j} \varphi ) (s,x/N) 
\big\{ \eta_s^M (x) - 1/2 - N^{-b} \varphi (s,x/N) \big\}^2 \Big] \\
\!\!\!\!\!\!\!\!\!\!\!\!\!\!\! && 
+\; \frac {4 N^{b}} {N^d}
\int_0^t ds \, \bb E_{\mu^N} \Big[ \sum_{j=1}^d 
\sum_{x \in \bb T_N^d} \gamma_j^M (\partial_{u_j} \varphi^2 ) (s,x/N) 
\big\{ \eta_s^M (x) - 1/2 - N^{-b} \varphi (s,x/N) \big\} \Big] \\
\!\!\!\!\!\!\!\!\!\!\!\!\!\!\! && 
+\; \frac {4  } {N^d}
\int_0^t ds \, \sum_{j=1}^d \sum_{x \in \bb T_N^d} 
\gamma_j^M (\partial_{u_j} \varphi ) (s,x/N) \varphi^2 (s,x/N) \;.
\end{eqnarray*}
As $N\uparrow\infty$, for each fixed $j$, $s$, the last term of this
expression converges to $4 \int_{\bb T^d} du \, \gamma_j \,
(\partial_{u_j} \varphi ) (s,u) \, \varphi^2 (s,u) = 0$. By Lemma
\ref{h09} below, the first term is bounded by $\gamma_0 \int_0^t ds
H_N(s) + o_N(1)$ for some finite constant $\gamma_0$. In the second
term, since $\varepsilon_N N^b$ vanishes as $N\uparrow\infty$ and
since, by \eqref{ass}, $N^b \ll M$, we may replace $(1 +
\varepsilon_N) \gamma_j^M$ by $\gamma_j$. The resulting expression
cancels with the third term of \eqref{hf06} because $\varphi$ is the
solution of \eqref{f11}. This proves Proposition \ref{h01} and
therefore Theorem \ref{mt3}.

We conclude this section with an estimate on the variance of the
density in terms of the relative entropy.

\begin{lemma}
\label{h09}
There exists $\gamma_0>0$ such that
\begin{equation*}
\bb E_{\mu^N} \Big[ \frac {N^{2b}} {N^d} 
\sum_{x \in \bb T_N^d} \Big\{ \eta_t^M (x) - (1/2)
- N^{-b} \varphi (t,x/N) \Big\}^2 \Big] \;\le\;
\gamma_0 H_N(t) \;+\; o_N(1)
\end{equation*}
for $0\le t\le T$. 
\end{lemma}

\begin{proof}
By the entropy inequality the expectation appearing in the statement
of the lemma is bounded above by
\begin{equation*}
\frac 1\gamma H_N(t) \;+\; \frac{N^{2b}}{\gamma N^d} \log
\bb E_{\nu_t^N} \Big[ \exp\Big\{ \gamma \sum_{x \in \bb T_N^d}
\Big( \eta_t^M (x) - (1/2) - N^{-b} \varphi (t,x/N) \Big)^2 
\Big\} \Big]
\end{equation*}
for every $\gamma >0$. By H\"older inequality, the second term is less
than or equal to
\begin{equation*}
\frac{N^{2b}}{\gamma N^d M^d} \sum_{x \in \bb T_N^d} \log
\bb E_{\nu_t^N} \Big[ \exp\Big\{ \gamma |\Lambda_M|
\Big( \eta_t^M (x) - (1/2) - N^{-b} \varphi (t,x/N) \Big)^2 \Big\}
\Big]\;. 
\end{equation*}
The above expectation is bounded uniformly in $N$ provided $\gamma$ is
small enough. The expression is thus bounded by $\gamma^{-1} N^{2b}
M^{-d}$, which concludes the proof of the lemma.
\end{proof}

\section{Proof of the Incompressible limit.}
\label{sec4}

Fix the reference measure $\nu^N_* = \nu^N_{(a_0,\bs 0)}$. Consider a
sequence of probability measures $\{\mu^N: N\ge 1\}$ satisfying the
assumptions of Theorem \ref{gs01}. A straightforward argument, similar
to the one which led to \eqref{f02}, shows that
\begin{equation*}
N^{2b-d} H_N(\mu^N | \nu^N_*) \;\le\; C_0
\end{equation*}
for some finite constant depending only on $a_0, \varphi_0, \bs
\varphi$. 

Denote by $f^N_t$ the Radon-Nikodym derivative $d\mu^N S_t^N /
d\nu^N_*$ and recall that $f^N_t$ solves the equation
\begin{equation*}
\partial_t f^N_t \;=\; \mc L_N^* f^N_t\;,
\end{equation*}
where $\mc L_N^*$ stands for the adjoint of $\mc L_N$ in
$L^2(\nu^N_*)$. By the previous estimate on the relative entropy of
$\mu^N$ with respect to $\nu^N_*$, we get that
\begin{equation}
\label{g16}
H_N(\mu^N S_t^N | \nu^N_*) \; +\; \int_0^t ds \, D_N(f^N_s)
\;\le\; C_0 N^{d-2b}\;,
\end{equation}
where $D_N$ stands for the Dirichlet form: $D_N(f) = \<f^{1/2}, (-\mc
L_N) f^{1/2}\>_{\nu^N_*}$.

Let $\psi_t^N = d \nu^N_t/ d\nu^N_*$. It follows from the explicit
formulas for the product measures $\nu^N_{t}$ that
\begin{equation*}
\log \psi_t^N = \sum_{x\in \bb T_N^d} \bs \lambda(t,x)
\cdot \bs I(\eta_x)\; -\; \sum_{x\in \bb T_N^d}
\log \frac {Z \big ( \bs \lambda(t,x)\big )}{Z(\bs 0) )}\;,
\end{equation*}
where $\bs \lambda(t,x):= \bs\Lambda( a_0 + N^{-b}\varphi_0, N^{-b}\bs
\varphi(t,x/N))$ and
\begin{equation*}
Z (\bs \lambda) \;=\; \sum_{\xi \in \{0,1\}^{\mc V}} 
\exp\big\{\bs \lambda \cdot \bs I(\xi) \big\}\;.
\end{equation*}

Let $H_N(t) = N^{2b - d} H_N(\mu^N S_t^N | \nu^N_t)$. With the
notation just introduced, we have that
\begin{equation*}
H_N(t) \;=\; N^{2b - d} \int f_t^N \log \frac{f_t^N}{\psi_t^N} 
\, d\nu^N_*\;.
\end{equation*}
Theorem \ref{gs01} follows from Gronwall lemma and the following
estimate.

\begin{proposition}
\label{gs06}
Fix a sequence of measures $\{\mu^N : N\ge 1\}$ satisfying the
assumptions of Theorem \ref{gs01}.  There exists $\gamma >0$ such that
\begin{equation*}
H_N(t) \;\le\; \gamma \int_0^t H_N(s) \, ds \;+\; o_N(1)
\end{equation*}
for all $0\le t\le T$.
\end{proposition}

The proof of Proposition \ref{gs06} is divided in several steps. We
begin with a well known upper bound for the entropy production.
\begin{equation}
\label{g14}
\frac{d}{dt} H_N(t) \;\le\; N^{2b-d}
\int f_t^N \frac{( \mc L_N^* - \partial_t) \psi_t^N}{\psi_t^N} 
\, d\nu^N_*\;.
\end{equation}

Next result is needed in to discard irrelevant terms on the right hand
side of the previous expression.

\begin{lemma}
\label{gs02}
Let $G : \bb T^d \to \bb R$ a continuous function and $\{\mu^N : N\ge
1\}$ a sequence of measures satisfying the assumptions of Theorem
\ref{gs01}. Then,
\begin{equation*}
\bb E_{\mu^N} \Big[ N^{-d} \sum_{x\in\bb
  T_N^d} G(x/N) I_k(\eta_x(t)) \Big] \;\le \; H_N(t) \;+\;
O(N^{-b})
\end{equation*}
for $0\le t\le T$, $1\le k\le d$. The lemma remains in force for $k=0$
if we replace $I_k(\eta_x(t))$ by $I_0(\eta_x(t)) - a_0$.
\end{lemma}

\begin{proof}
Fix $1\le k\le d$. We may replace $I_k(\eta_x(t))$ by $I_k(\eta_x(t))
- N^{-b} \varphi_k(t,x/N)$ paying a price of order $N^{-b}$. It
remains to apply the entropy inequality with respect to measure
$\nu_t$, which is product, and perform a second order Taylor
expansion. 
\end{proof}

A long and tedious computation gives that
$(\psi_t^N)^{-1} ( \mc L_N^* - \partial_t) \psi_t^N$ is equal to
\begin{eqnarray}
\label{g10}
\!\!\!\!\!\!\!\!\!\!\!\!\!\!\! && \nonumber
\frac{4}{B N^b} \sum_{k=1}^d \sum_{x \in \bb T_N^d} (\Delta
\varphi_k)(t,x/N) I_k(\eta_x) \;+\; \frac 4B \sum_{j,k=1}^d
\sum_{x \in \bb T_N^d}
(\partial_{u_j} \varphi_k) (t,x/N) \tau_x W_{k,j}^{*,M} \\
\!\!\!\!\!\!\!\!\!\!\!\!\!\!\! && 
+\; \frac{16}{B^2N^{2b}} \sum_{i,j=1}^d \sum_{k,\ell =1}^d  
\sum_{x \in \bb T_N^d} (\partial_{u_i} \varphi_k) (t,x/N)
(\partial_{u_j} \varphi_\ell) (t,x/N) \tau_x V_{i,j}^{k,\ell,M} \\
\!\!\!\!\!\!\!\!\!\!\!\!\!\!\! && \nonumber
-\; \frac 4{BN^{b}} \sum_{k=1}^d (\partial_t \varphi_k) (t,x/N)
\Big\{ I_k(\eta_x) - \frac{\varphi_k (t,x/N)}{N^b} \Big\}\; 
+ \; R_N(t) \;+\; o(N^{d-2b})\; .
\end{eqnarray}
In this formula, $W_{k,j}^{*,M}$ and $V_{i,j}^{k,\ell,M}$ stand for
\begin{eqnarray*}
\!\!\!\!\!\!\!\!\!\!\!\!\!\!\! && 
W_{k,j}^{*,M} \;=\; \frac {A_M}{M^{d+1}} \sum_{v\in\mc V} v_k\, 
\sum_{z\in\Lambda_M} q_M(-z, v) \, z_j \, \eta(0,v)\, \{ 1- \eta(z,v)\} \\
\!\!\!\!\!\!\!\!\!\!\!\!\!\!\! && 
\quad 
V_{i,j}^{k,\ell,M} \;=\; \frac {A_M}{M^{d+2}} \sum_{v\in\mc V} 
v_k \, v_\ell\,  \sum_{z\in \Lambda_M} z_i \, z_j \,
\eta(0,v) \, \{1-\eta(z,v)\}\;. 
\end{eqnarray*}
Since the density $\psi_t^N$ is a function of the conserved quantities
$\bs I$, the collision part of the generator is irrelevant in the
previous computation.  We used repeatedly Lemma \ref{gs02} and the
fact that $b<a$, which follows from \eqref{g09}, to discard
superfluous terms. The remainder $o(N^{d-2b})$ should be understood as
an expression whose expectation with respect to $\mu^N S_t^N$
integrated in time is of order $o(N^{d-2b})$, while $R_N(t)$ is an
expression which multiplied by $N^{2b-d}$ is bounded by $H_N(t) +
O(N^{-b})$ in virtue of Lemma \ref{gs02}.

If we replace $I_k(\eta_x)$ in the first term of \eqref{g10} by
$I_k(\eta_x) - \varphi_k (t,x/N) N^{-b}$ and $V_{i,j}^{k,\ell,M}$ by
$V_{i,j}^{k,\ell,M} - (B/4) \delta_{i,j} \delta_{k,\ell}$, as
$N\uparrow\infty$, the expressions added when multiplied by $(B/4)
N^{2b-d}$ converge to 
\begin{equation*}
\sum_{k=1}^d \int_{\bb T^d} (\Delta \varphi_k) \, \varphi_k \;+\;
\sum_{j,k=1}^d \int_{\bb T^d} (\partial_{u_j} \varphi_k)^2 \;=\; 0\;.
\end{equation*}

Therefore, in view of \eqref{g14}, \eqref{g10} and Lemma \ref{gs02},
the time derivative of the renormalized entropy $H_N(t)$ is bounded
above by
\begin{eqnarray}
\label{g11}
\!\!\!\!\!\!\!\!\!\!\!\!\!\!\! && H_N(t)+
\bb E_{\mu^N} \Big[
\frac {4 N^{2b}} {BN^d}  \sum_{j,k=1}^d \sum_{x \in \bb T_N^d}
(\partial_{u_j} \varphi_k) (t,x/N) \tau_x W_{k,j}^{*,M} (t) \Big] \\
\!\!\!\!\!\!\!\!\!\!\!\!\!\!\! && \nonumber
+\; \bb E_{\mu^N} \Big[ \frac{16}{B^2N^{d}} \sum_{i,j=1}^d \sum_{k,\ell =1}^d  
\sum_{x \in \bb T_N^d} (\partial_{u_i} \varphi_k) (t,x/N)
(\partial_{u_j} \varphi_\ell) (t,x/N) 
\tau_x  \hat V_{i,j}^{k,\ell,M} (t) \Big] \\
\!\!\!\!\!\!\!\!\!\!\!\!\!\!\! && \nonumber
+\; \bb E_{\mu^N} \Big[ 
\frac {4 N^b}{BN^{d}} \sum_{k=1}^d \sum_{x\in\bb T^d}
( \Delta \varphi_k - \partial_t \varphi_k) (t,x/N)
\Big\{ I_k(\eta_x(t)) - \frac{\varphi_k (t,x/N)}{N^b} \Big\}\Big] \; 
+ \; o_N(1)\; ,
\end{eqnarray}
where $\hat V_{i,j}^{k,\ell,M} = V_{i,j}^{k,\ell,M} - (B/4)
\delta_{i,j} \delta_{k,\ell}$.

We now use the ergodicity to replace the functions $W_{k,j}^{*,M}$ and
$\hat V_{i,j}^{k,\ell,M}$ by their projections on the conserved
quantities. For $s\ge 0$ and $x$ in $\bb Z^d$, denote by $\bs
I^M(s,x)$ the average at time $s$ of the conserved quantities over a
cube $\Lambda_M$ centered at $x$:
\begin{equation*}
\bs I^M(s,x) \;=\; \frac 1{|\Lambda_M|} \sum_{y\in x + \Lambda_M}
\bs I (\eta_y (s))\;.
\end{equation*}
To keep notation simple, let $\bs I^M_s:= \bs I^M(s,0)$. 

Recall the definition of the canonical measures $\nu_{\Lambda_M,\bs
  i}$ presented in Subsection \ref{sec3.3}. Since we assumed in
\eqref{g12} the global dynamics restricted to a cube of length $M$ to
have a spectral gap of order $M^\kappa$ and since $a+b > 1 -
[2/(d+\kappa)]$, $a + (\kappa -2 /\kappa)b > 1 - (2/\kappa)$,
repeating the arguments presented in the proof of Lemma \ref{s1} and
taking advantage of the estimate \eqref{g16} we derive the so-called
one block estimate. In this lemma, the collision part of the dynamics,
also speeded up by $N^2$, plays an important role.

\begin{lemma}
\label{gs04}
For every $t\ge 0$, every $1\le j,k\le d$ and every continuous
function $G: [0,T]\times \bb T^d\to\bb R$,
\begin{equation*}
\limsup_{N\to\infty} \bb E_{\mu^N} \Big[
\,\Big\vert \int_0^t ds\, \frac {N^{2b}}{N^d} \sum_{x\in\bb T_N^d}
G(s, x/N) \tau_x \Big\{ W_{k,j}^{*,M}(s)
-E_{\nu_{\Lambda_M,\bs I^M_s}}[W_{k,j}^{*,M}]\Big\} 
\Big\vert \, \Big] \;=\; 0\;.
\end{equation*}
\end{lemma}

Since $\nu_{\Lambda_M,\bs i}$ is the counting measure,
\begin{equation*}
E_{\nu_{\Lambda_M,\bs i}} \Big[ W_{k,j}^{*,M} \Big]
\;=\; - \sum_{v\in\mc V} v_k v_j E_{\nu_{\Lambda_M,\bs i}} 
\Big[ \eta(0,v) [1-\eta(e_1,v)] \Big]
\end{equation*}
because $A_M M^{-(d+1)} \sum_{z\in\Lambda_M} q_M(z,v) z_j = v_j$.  In
the previous formula, site $e_1$ can be replaced by any site of
$\Lambda_M$ different from the origin. Since $N^{2b} \ll M^d$, by the
equivalence of ensembles, stated in Proposition \ref{hs1} below, we
can replace the expectation with respect to the canonical measure by
the expectation with respect to the grand canonical measure paying a
price of order $o_N(1)$.

For $1\le j, k\le d$, let 
\begin{equation*}
R_{j,k}(\rho,\bs p) \; := \; 
E_{\mu^N_{\Lambda(\rho,\bs p)}} \big[ W_{k,j}^{*,M} \big]
\; = \; - \sum_{v\in\mc V} v_k\, v_j\,  
\chi (\theta_v(\bs \Lambda(\rho,\bs p))) \;,
\end{equation*}
where $\theta_v(\cdot)$ is defined in \eqref{hf2}.  Up to this point,
we replaced the first expectation in \eqref{g11} by
\begin{equation*} 
\bb E_{\mu^N} \Big[ \frac {4 N^{2b}} {B N^d}  \sum_{j,k=1}^d 
\sum_{x \in \bb T_N^d} (\partial_{u_j} \varphi_k) (t,x/N) 
R_{j,k} (\bs I^M (t,x))\Big]\; +\; o_N(1)\;.
\end{equation*}

Since $\bs \Lambda(a_0, \bs 0) = \bs 0$ and $\theta_v(\bs 0) = 1/2$,
$R_{j,k}(a_0,\bs 0)= -(B/4) \delta_{j,k}$. On the other hand, since
$\varphi$ is divergence free,
\begin{equation*}
\frac {4 N^{2b}} {B N^d}  \sum_{j,k=1}^d \sum_{x \in \bb T_N^d}
(\partial_{u_j} \varphi_k) (t,x/N) R_{j,k} (a_0,\bs 0)
\;=\; \frac {- N^{2b}} {N^d}  \sum_{j=1}^d \sum_{x \in \bb T_N^d}
(\partial_{u_j} \varphi_j) (t,x/N)
\end{equation*}
vanishes for each fixed $N$. We may therefore add this expression to
the previous expectation to obtain that the first term in \eqref{g11} is
equal to
\begin{equation}
\label{g15}
\bb E_{\mu^N} \Big[ \frac {4 N^{2b}} {B N^d}  \sum_{j,k=1}^d 
\sum_{x \in \bb T_N^d} (\partial_{u_j} \varphi_k) (t,x/N) 
\big\{ R_{j,k} (\bs I^M (t,x)) - R_{j,k} (a_0,\bs 0) \big\}\Big]
\; +\; o_N(1)\;.
\end{equation}

The same arguments show that we can replace $\hat V_{i,j}^{k,\ell,M}$
in the second term of \eqref{g11} by its expectation with respect to
the grand canonical measure. The proof is even simpler due to the
absence of the factor $N^{2b}$ in front of the sum.  Since
$R_{j,k}(a_0,\bs 0)= -(1/4) \delta_{j,k} B$,
\begin{eqnarray*}
E_{\mu^N_{\bs \Lambda(\rho,\bs p)}} \big[ \hat V_{i,j}^{k,\ell,M} \big]
\;=\; - \delta_{i,j}
\big[ R_{k,\ell}(\rho,\bs p)-R_{k,\ell}(a_0,\bs 0) \big]\; . 
\end{eqnarray*}
The one-block estimate permits therefore to replace the second
expectation in \eqref{g11} by
\begin{equation}
\label{V}
-\;\bb E_{\mu^N} \Big[\frac{16}{B^2N^{d}}\sum_{i,k,\ell =1}^d  
\sum_{x \in \bb T_N^d} (\partial_{u_i} \varphi_{k,\ell}) (t,x/N) 
\big\{ R_{k,\ell}(\bs I^M(t,x))-R_{k,\ell}(a_0,\bs 0) \big\} 
\Big]\;,
\end{equation}
where $(\partial_{u_i} \varphi_{k,\ell}) (t,x/N)=(\partial_{u_i}
\varphi_k)(t,x/N) (\partial_{u_i} \varphi_\ell) (t,x/N)$.

It is now clear that \eqref{V} is a term of lower order than
\eqref{g15}. We therefore only need to estimate the latter. Fix an
arbitrary $\epsilon >0$. Since $R_{j,k}$ is a bounded function, the
integral in \eqref{g15} when restricted to $|\bs I^M (t,x)- (a_0,\bs
0)|>\epsilon $ is bounded above by
 
\begin{equation}
\label{h15}
\bb E_{\mu^N} \Big[ \frac {C_0 N^{2b}} {\epsilon^3 N^d}
\sum_{x \in \bb T_N^d} 
\big|\bs I^M (t,x) - (a_0,\bs 0) \big|^3 \Big]\;,
\end{equation}
where $C_0$ is a constant depending on $\mc V$ and $\bs \varphi$. In
the expression above we may replace $(a_0,\bs 0)$ by
$(a_0+N^{-b}\varphi_0, N^{-b}\bs \varphi(t,x/N))$ paying a price of
order $N^{-b}$. Since $I^M (t,x)$ belongs to a compact set the
expression obtained after replacing is bounded above by
\begin{equation*}
\bb E_{\mu^N} \Big[ \frac {C(\mc V)C_0 N^{2b}}{\epsilon^3 N^d}
\sum_{x \in \bb T_N^d}
\Big| \bs I^M (t,x) - (a_0+\frac{\varphi_0}{N^b},\frac{\bs
  \varphi(t,x/N)}{N^b}) \Big|^2 \Big]\;. 
\end{equation*}
By Lemma \ref{gs09} below, this expression is bounded by $\gamma_0
H_N(t) + o_N(1)$ for some $\gamma_0>0$.

In order to deal with the integral \eqref{g15} on $|\bs
I^M(t,x)-(a_0,\bs 0)|\le \epsilon$ we perform a Taylor expansion of
$R_{j,k}$. The first term in the expansion vanishes because the
gradient of $R_{j,k}$ vanishes at $(a_0,\bs 0)$. The contribution of
the second order terms is
\begin{equation*}
\frac {4 N^{2b}} {B N^d}  \sum_{j,k=1}^d \sum_{x \in \bb T_N^d}
(\partial_{u_j} \varphi_k) (t,x/N) \sum_{v\in\mc V} v_k\, 
v_j\, \Big\{  \frac{I_0^L - a_0}{|\mc V|}  \;+\; \frac 1B
\sum_{\ell=1}^d  v_\ell\, I^L_\ell  \Big\}^2 \;.
\end{equation*}
Expanding the square, the term in $I_0^L$ vanishes because $\bs \varphi$
is divergence free and the cross product vanishes because $\mc V$ is
symmetric. This sum is therefore equal to
\begin{equation*}
\frac {4 N^{2b}} {B^3 N^d}  \sum_{j,k,\ell, m =1}^d \sum_{x \in \bb T_N^d}
(\partial_{u_j} \varphi_k) (t,x/N) \sum_{v\in\mc V} v_k\, 
v_j\,  v_\ell\,  v_m\, I^L_\ell \, I^L_m  \;.
\end{equation*}

Replacing $I_\ell^L$ by $I_\ell^L - \varphi_\ell(t,x/N)
N^{-b}$, we may rewrite the previous expression as the sum of three
kind of terms. The first one, the $0$ order term in $I$, consists
simply in replacing $I^L_\ell$ by $\varphi_\ell N^{-b}$. As $N$ tends
to infinity, this term converges to
\begin{equation*}
\frac{4(D-3C)}{B^3} \sum_{j=1}^d \int_{\bb T^d} (\partial_{u_j}
\varphi_j) \varphi_j^2 \;+\; \frac {4C}{B^3} \sum_{j,k=1}^d 
\int_{\bb T^d} (\partial_{u_j} \varphi_k^2) \varphi_j\;.
\end{equation*}
An integration by parts shows that this expression vanishes because
$\bs \varphi$ is divergence free. The linear term in $I$ cancels with
the last term of \eqref{g11} because $\bs \varphi$ is the solution of the
Navier-Stokes equation \eqref{g05}. Remains the quadratic term in $I$,
equal to
\begin{eqnarray*}
\!\!\!\!\!\!\!\!\!\!\!\!\!\! &&
\frac {4 (D-3C) N^{2b}} {B^3 N^d} \sum_{j=1}^d \sum_{x \in \bb T_N^d}
(\partial_{u_j} \varphi_j) (t,x/N) \Big\{ I^L_j (t,x) - 
\frac{\varphi_j (t,x/N)}{N^b} \Big\}^2 \; + \\
\!\!\!\!\!\!\!\!\!\!\!\!\!\! &&
\frac {8 C N^{2b}} {B^3 N^d} \sum_{j,k=1}^d \sum_{x \in \bb T_N^d}
(\partial_{u_j} \varphi_k) (t,x/N) \Big\{ I^L_j (t,x) - 
\frac{\varphi_j (t,x/N)}{N^b} \Big\} \Big\{ I^L_k (t,x) - 
\frac{\varphi_k (t,x/N)}{N^b} \Big\}
\end{eqnarray*}
because $\bs \varphi$ is divergence free. By Lemma \ref{gs09} below,
this expression is bounded by $\gamma_0 H_N(t) + o_N(1)$ for some
$\gamma_0>0$.

Finally, we consider the remainder in the Taylor expansion. Since
$R_{j,k}$ is smooth, we can choose $\epsilon$ small enough for the
third derivative of $R_{j,k}$ to be bounded in an
$\epsilon$-neighborhood of $(a_0, \bs 0)$ by a finite constant $C_0$
depending on $\mc V$ and $\bs\varphi$. In particular, the remainder is
bounded above by
\begin{equation*}
\frac {C_0 N^{2b}} {N^d} \sum_{x \in \bb T_N^d}
\big|\bs I^M (t,x) - (a_0, \bs 0)\big|^3\;.
\end{equation*}
The same arguments used to estimate \eqref{h15} prove that this
expression is bounded by $\gamma_0 H_N(t) + o_N(1)$ for some
$\gamma_0>0$. This concludes the proof of Proposition \ref{gs06}.
\smallskip

We conclude the section with an estimate repeatedly used in the proof
of Proposition \ref{gs06}. We need here again the assumption that
$N^{2b}\ll M^d$.

\begin{lemma}
\label{gs09}
There exists $\gamma_0>0$ such that
\begin{equation*}
\bb E_{\mu^N} \Big[ \frac {N^{2b}} {N^d} 
\sum_{x \in \bb T_N^d} \Big\{ I^M_j (t,x) - 
\frac{\varphi_j (t,x/N)}{N^b}\Big\}^2 \Big] \;\le\;
\gamma_0 H_N(t) \;+\; o_N(1)
\end{equation*}
for $1\le j\le d$, $0\le t\le T$. The statement remains in force for
$j=0$ if $\varphi_j (t,x/N) N^{-b}$ is replaced by $a_0 + \varphi_0
N^{-b}$.
\end{lemma}

\begin{proof}
By the entropy inequality the expectation appearing in the statement
of the lemma is bounded above by
\begin{equation*}
\frac 1\gamma H_N(t) \;+\; \frac{N^{2b}}{\gamma N^d} \log
\bb E_{\nu_t^N} \Big[ \exp\Big\{ \gamma \sum_{x \in \bb T_N^d}
\Big( I^M_j (t,x) - \frac{\varphi_j (t,x/N)}{N^b}\Big)^2 \Big\} \Big]
\end{equation*}
for every $\gamma >0$. By H\"older inequality, the second term is less
than or equal to
\begin{equation*}
\frac{N^{2b}}{\gamma N^d M^d} \sum_{x \in \bb T_N^d} \log
\bb E_{\nu_t^N} \Big[ \exp\Big\{ \gamma |\Lambda_M|
\Big\{ I^M_j (t,x) - \frac{\varphi_j (t,x/N)}{N^b}\Big\}^2 \Big\}
\Big]\;. 
\end{equation*}
The above expectation is bounded uniformly in $N$ provided $\gamma$ is
small enough. The expression is thus bounded by $\gamma^{-1} N^{2b}
M^{-d}$, which concludes the proof of the lemma.
\end{proof}

\section{Spectral gap for stochastic lattice gases}
\label{gap}

We prove in this section a spectral gap of polynomial order for the
generator of the stochastic lattice gas. We consider a slightly
different process, in which the exclusion dynamics allows particles to
jump to any site of $\Lambda_M$ at rate $M^{-(d+2)}$. We do not
require, therefore, the jump to be of size smaller than $M$. Of
course, the Dirichlet forms of both dynamics are equivalent and the
result stated in Proposition \ref{gs10} extends to the original
dynamics.

Fix $M\ge 1$ and consider the process restricted to the cube
$\Lambda_M$ without the factor $N^2$. The generator of the process,
denoted by $\mc L_M$, can be written as $\mc L_M^{ex} + \mc L_M^c$,
where
\begin{eqnarray*}
\!\!\!\!\!\!\!\!\!\!\!\!\! &&
(\mc L_M^{ex} f) (\eta) \;=\; \frac 1{M^{d+2}}
\sum_{v \in \mc V} \sum_{x,z \in \Lambda_M}
\eta(x,v) \, [1-\eta(z,v)] \, [f(\eta^{x, z, v}) -
f(\eta)]\;, \\ 
\!\!\!\!\!\!\!\!\!\!\!\!\! && \qquad
\mc L_M^c f(\eta) \;=\; \sum_{y \in \Lambda_M}
\sum_{q \in \mc Q} p(y,q,\eta) \, [f(\eta^{y,q}) - f(\eta)]
\end{eqnarray*}
and $p(y,q,\eta)$ is defined at the beginning of Section \ref{sec1}

For each fixed $\bs i$ in $\mf V_M$, recall that we denote by
$\nu_{\Lambda_M, \bs i}$ the invariant measure concentrated on
configurations $\eta$ of $(\{0,1\}^{\mc V})^{\Lambda_M}$ such that
$\bs I^M(\eta) = \bs i$. An elementary computations shows that
\begin{eqnarray}
\label{g24}
\!\!\!\!\!\!\!\!\!\!\!\!\!\!\!\! &&
\< f, - \mc L_M^{ex} f\>_{\nu_{\Lambda_M,\bs i}} \;=\;
\frac {1}{4 M^{d+2}} \sum_{v\in\mc V}
\sum_{x,y\in \Lambda_M} E_{\nu_{\Lambda_M,\bs i}} 
\Big[ \{f(\xi^{x, y, v}) - f(\xi)\} ^2 \Big] \; , \\
\!\!\!\!\!\!\!\!\!\!\!\!\!\!\!\! && \quad
\< f, - \mc L_M^c f\>_{\nu_{\Lambda_M,\bs i}}
\;=\;  \frac 12 \sum_{q\in\mc Q} \sum_{x\in \Lambda_M}
E_{\nu_{\Lambda_M,\bs i}} \Big[ p(x, q, \xi) 
\{f(\xi^{x, q}) - f(\xi)\} ^2 \Big] 
\nonumber
\end{eqnarray}

Denote by $E_{\nu}[f;f]$ the variance of $f$ with
respect to a measure $\nu$ and by $\< \cdot , \cdot \>_{\nu}$ the
inner product in $L^2(\nu)$.

\begin{proposition}
\label{gs10}
There exists a finite constant $C_2$, depending only on $\mc V$, $\mc
Q$, such that
\begin{equation*}
E_{\nu_{\Lambda_M, \bs i}} [ f;f] \;\le\; C_2 M^{2+ 3d + 2 d^2} 
\< f, - \mc L_M f\>_{\nu_{\Lambda_M,\bs i}}
\end{equation*}
for all $f$ in $L^2(\nu_{\Lambda_M,\bs i})$, all
$\bs i$ in $\mf V_M$ and all $M\ge 1$. 
\end{proposition}

The proof of this proposition relies on estimates on the Dirichlet
forms associated to $\mc L_M^{ex}$ and $\mc L_M^c$.  Denote by
$\tilde{\mc L}^c_M$ the generator of a dynamics in which collisions
between particles at different sites are allowed:
\begin{equation*}
(\tilde{\mc L}^c_M f)(\eta) \;=\; \frac 1{|\Lambda_M|^3} 
\sum_{x_1, \dots, x_4  \in \Lambda_M}
\sum_{q \in \mc Q} p(\bs x, q, \eta) \, \{f(\eta^{\bs x,q}) - f(\eta)\}\;,
\end{equation*}
where, for $q=(u,v,u',v')$, $\bs x = (x_1, \dots, x_4)$,
\begin{equation*}
p(\bs x, q, \eta) \;=\; \eta_{x_1}(u) \, \eta_{x_2}(v) 
\, [1-\eta_{x_3}(u')] \, [1-\eta_{x_4}(v')]
\end{equation*}
and where $\eta^{\bs x, q}$ is the configuration $\eta$ in which the
occupation variables $\eta_{x_1}(u)$, $\eta_{x_2}(v)$,
$\eta_{x_3}(u')$, $\eta_{x_4}(v')$ are flipped.

The first lemma of this section states that the Dirichlet forms
associated to $\mc L^c_M$ and to $\tilde{\mc L}^c_M$ are comparable
and that the Dirichlet form of the conditional expectation of a
function with respect to the total number of particles with fixed
velocity can be estimated by the Dirichlet form of the original
function. For each $v$ in $\mc V$, let $K_v$ be the total number of
particles with velocity $v$ in $\Lambda_M$:
\begin{equation*}
K_v \;=\; K_v(\eta) \;=\; \sum_{x\in \Lambda_M} \eta(x,v)\;.
\end{equation*}

\begin{lemma}
\label{gs12}
There exists a finite constant $C_2$, depending only on
$\mc V$, such that
\begin{equation*}
\< - \tilde{\mc L}^c_M f, f\>_{\nu_{\Lambda_M,\bs i}} \;\le\;
C_2 M^2\, \< - \mc L^{ex}_M f, f\>_{\nu_{\Lambda_M,\bs i}} 
\;+\; C_2 \< - \mc L^c_M f, f\>_{\nu_{\Lambda_M,\bs i}} 
\end{equation*}
for every $f$ in $L^2(\nu_{\Lambda_M,\bs i})$. Moreover, let
\begin{equation*}
F \;=\; E_{\nu_{\Lambda_M,\bs i}}  
\big[ f\, \big\vert\, \{K_v, v\in\mc V\} \big]\;.
\end{equation*}
Then,
\begin{equation*}
\< - \mc L^c_M F, F\>_{\nu_{\Lambda_M,\bs i}} \;\le\;
\< - \tilde{\mc L}^c_M f, f\>_{\nu_{\Lambda_M,\bs i}}
\end{equation*}
for every $f$ in $L^2(\nu_{\Lambda_M,\bs i})$.
\end{lemma}

\begin{proof}
An elementary computation shows that 
\begin{equation*}
\< - \tilde{\mc L}^c_M f, f\>_{\nu_{\Lambda_M,\bs i}} \;=\;
\frac 1{2 |\Lambda_M|^3} \sum_{q\in \mc Q}
\sum_{x_1, \dots, x_4 \in\Lambda_M} 
E_{\nu_{\Lambda_M,\bs i}} \Big[ p(\bs x, q, \xi) 
\{f(\xi^{\bs x, q}) - f(\xi)\} ^2 \Big]\;.
\end{equation*}
for $f$ in $L^2(\nu_{\Lambda_M,\bs i})$. Fix $q$ and $\bs x$. We construct a
path from $\xi$ to $\xi^{\bs x, q}$ with jumps and collisions of
particles in the same site in the following way. Assume that the set
$\mc V$ has been ordered: $\mc V = \{v_1, \dots, v_n\}$ and, without
loss of generality, that $q=(v_1, \dots, v_4)$. We first exchange the
occupation variable $\xi_{x_2}(v_2)$, $\xi_{x_1}(v_2)$; than
$\xi_{x_3}(v_3)$, $\xi_{x_1}(v_3)$ and finally $\xi_{x_4}(v_4)$,
$\xi_{x_1}(v_4)$. At this point we may perform the collision at site
$x_1$ and move back the particles and holes to their final positions
in the reversed order.

The total length of the path is at most $7$. Denote by $\zeta_0 = \xi,
\dots, \zeta_\ell = \xi^{\bs x, q}$ the successive
configurations. Writing $\{f(\xi^{\bs x, q}) - f(\xi)\}$ as
$\sum_j \{f(\zeta_{j+1}) - f(\zeta_j)\}$, applying Schwarz inequality,
reversing the order of the summations and estimating the total number
of configurations whose path jumps from $\zeta_j$ to $\zeta_{j+1}$, we
obtain that for each $q=(v_1, \dots, v_4)$,
\begin{eqnarray*}
\!\!\!\!\!\!\!\!\!\!\!\!\! &&
\sum_{x_1, \dots, x_4 \in\Lambda_M} 
E_{\nu_{\Lambda_M, \bs i}} \Big[ p(\bs x, q, \xi) 
\{f(\xi^{\bs x, q}) - f(\xi)\} ^2 \Big]  \\
\!\!\!\!\!\!\!\!\!\!\!\!\! && \qquad\qquad\qquad
\;\le\; C_0 |\Lambda_M|^2 \sum_{i=1}^4 \sum_{x,y\in \Lambda_M}
E_{\nu_{\Lambda_M,\bs i}} \Big[ \{f(\xi^{x, y, v_i}) - f(\xi)\} ^2 \Big] \\
\!\!\!\!\!\!\!\!\!\!\!\!\! && \qquad\qquad\qquad
\;+\; C_0 |\Lambda_M|^3 \sum_{x\in \Lambda_M}
E_{\nu_{\Lambda_M, \bs i}} \Big[ p(x, q, \xi) \{f(\xi^{x, q}) - f(\xi)\} ^2
\Big] \;.
\end{eqnarray*}
for some finite constant $C_0$. In particular, summing over $q$ in
$\mc Q$ and dividing by $2 |\Lambda_M|^3$, we obtain that
\begin{eqnarray*}
\< - \tilde{\mc L}^c_M f, f\>_{\nu_{\Lambda_M,\bs i}} &\le&
\frac {C_2}{|\Lambda_M|} \sum_{v\in\mc V}
\sum_{x,y\in \Lambda_M} E_{\nu_{\Lambda_M, \bs i}} 
\Big[ \{f(\xi^{x, y, v}) - f(\xi)\} ^2 \Big] \\
&+&  C_2 \sum_{q\in\mc Q} \sum_{x\in \Lambda_M}
E_{\nu_{\Lambda_M,\bs i}} \Big[ p(x, q, \xi) 
\{f(\xi^{x, q}) - f(\xi)\} ^2 \Big] 
\end{eqnarray*}
for some finite constant $C_2$ depending on $\mc V$. This
expression is bounded by $$C_2 \{ \< - \mc L^c_M f, f\>_{\nu_{\Lambda_M,\bs i}}
+ M^2 \< - \mc L^{ex}_M f, f\>_{\nu_{\Lambda_M,\bs i}} \}$$ in view of
\eqref{g24}. This concludes the proof of the first statement of the
lemma. We turn now to the second.

Fix $x$ in $\Lambda_M$ and $q=(u,v,u',v')$ in $\mc Q$. An elementary
computation shows that 
\begin{equation*}
F(\eta^{x,q}) - F(\eta) \;=\;
\frac 1{Z_q (\eta)} \sum_{x_1, \dots, x_4 \in\Lambda_M} 
E\Big[ p(\bs x, q, \xi) \{f(\xi^{\bs x, q}) - f(\xi)\}\, \big\vert\, 
\{K_v(\eta) : v\in \mc V\} \Big] \; ,
\end{equation*}
where 
\begin{eqnarray*}
Z_q (\eta) \;=\; K_u(\eta) K_v(\eta) (|\Lambda_M| - K_{u'})
(|\Lambda_M| - K_{v'})\;.
\end{eqnarray*}
In particular, by Schwarz inequality,
\begin{eqnarray*}
\!\!\!\!\!\!\!\!\!\!\!\!\!\! &&
E_{\nu_{\Lambda_M,\bs i}} \Big[ p(x,q,\eta) \{ F(\eta^{x,q}) - F(\eta) \}^2
\Big] \;\le\; \\
\!\!\!\!\!\!\!\!\!\!\!\!\!\! && \quad
\sum_{x_1, \dots, x_4 \in\Lambda_M} 
E_{\nu_{\Lambda_M,\bs i}} \Big[ \frac {p(x,q,\eta)} {Z_q (\eta)}  
E\Big[ p(\bs x, q, \xi) \{f(\xi^{\bs x, q}) - f(\xi)\} ^2 \, \big\vert\, 
\{K_v(\eta) : v\in \mc V\} \Big] \Big]\;.
\end{eqnarray*}
Taking conditional expectation with respect to $\{K_v(\eta) : v\in \mc
V\}$, summing over $x$ $q$ and dividing by $2$, the previous
expression becomes
\begin{equation*}
\frac 1{2 |\Lambda_M|^3} \sum_{q\in \mc Q}
\sum_{x_1, \dots, x_4 \in\Lambda_M} 
E_{\nu_{\Lambda_M,\bs i}} \Big[ p(\bs x, q, \xi) 
\{f(\xi^{\bs x, q}) - f(\xi)\} ^2 \Big]\;=\;
\< - \tilde{\mc L}^c_M f, f\>_{\nu_{\Lambda_M,\bs i}}\;.
\end{equation*}
This concludes the proof of the lemma.
\end{proof}

\noindent{\sl Proof of Proposition \ref{gs10}.}
Fix $M\ge 1$, $\bs i$ in $\mf V_M$ and a function $f$ in
$L^2(\nu_{\Lambda_M,\bs i})$.  Denote the conditional expectation of $f$ with
respect to $\bs K = \{K_v, v\in\mc V\}$ by $F(\bs K)$:
\begin{equation*}
F(\bs K) \;=\; E_{\nu_{\Lambda_M,\bs i}} 
\big[ f\, \big\vert\, \{K_v, v\in\mc V\} \big]\;.
\end{equation*}
By orthogonality,
\begin{equation}
\label{g17}
E_{\nu_{\Lambda_M,\bs i}} [ f;f] \;=\;
E_{\nu_{\Lambda_M,\bs i}} \big[ \{f - F(\bs K)\}^2 \big] \;+\;
E_{\nu_{\Lambda_M,\bs i}} [ F; F ]\;.
\end{equation}
Using only the exclusion part of the dynamics, since particles jump
uniformly over the cube $\Lambda_M$ with rate $M^{-(d+2)}$, by \cite{q},
\begin{equation*}
E_{\nu_{\Lambda_M,\bs i}} \big[ \{f - F(\bs K)\}^2 \big] \;\le\;
C_0 M^2 \< f, - \mc L^{ex}_M f\>_{\nu_{\Lambda_M,\bs i}}
\end{equation*}
for some finite universal constant $C_0$.

To estimate the second piece on the right hand side of \eqref{g17},
note that the exclusion part is irrelevant, while the Dirichlet form
associated to the collision part can be written as
\begin{eqnarray*}
\!\!\!\!\!\!\!\!\!\!\! &&
\< F, - \mc L^{c}_M F\>_{\nu_{\Lambda_M,\bs i}} \;=\; \\
\!\!\!\!\!\!\!\!\!\!\! && \qquad
|\Lambda_M|^{-3} \sum_{q\in \mc Q} \<\, K_v K_w (|\Lambda_M|- K_{v'})
(|\Lambda_M|- K_{v'}) \{F(\bs K^q)-F(\bs K)\}^2\,\>_{\nu_{\Lambda_M,\bs i}}\;.
\end{eqnarray*}
Denote by $\Omega_{M,\bs i}$ the state space of velocities
on $\Lambda_M$,
\begin{equation*}
\Omega_{M,\bs i} \;=\; \Big\{ \bs K=(K_1, \dots , K_v): 
\sum_v \bs K= |\Lambda_M|\bs i\Big\}\;,
\end{equation*}
and by $\bar \nu_{\Lambda_M,\bs i}$ the invariant state $\nu_{\Lambda_M,\bs i}$
projected on $\Omega_{M,\bs i}$. An elementary computation shows that
\begin{equation*}
\bar \nu_{\Lambda_M,\bs i} (\bs K) \;=\; \frac 1{Z_{M,\bs i}} 
\prod_{v\in\mc V}  
\begin{pmatrix}
  |\Lambda_M| \\ K_v
\end{pmatrix}
\end{equation*}
for some renormalizing constant $Z_{M,\bs i}$.

Consider from now on model I. Set $|\Lambda_M|\bs i:=(I_0,\dots,I_d)$.
Suppose without loss of generality that $K_{-e_j} \le K_{e_j}$ for
$1\le j\le d$ and let $K_j=K_{-e_j}$. Since $I_j = K_{e_j} -
K_{-e_j}$, $\bs K$ can be recovered from $\bs i$, $K_j$. We may
therefore ignore $(K_{e_1}, \dots, K_{e_d})$ and assume that
$(K_1,\dots, K_d)$ is evolving on the hyperplane
\begin{equation*}
\mc H \;=\; \mc H_{M,\bs i} 
\;=\; \Big\{(K_1, \dots , K_d) : K_j\ge 0, \, 2\sum_j K_j = 
I_0 - \sum_j I_j \Big\}\;.
\end{equation*}
On the set $\mc H$ the measure $\bar \nu_{\Lambda_M,\bs i}$ becomes
\begin{equation*}
\bar \nu_{\Lambda_M,\bs i} (K_1, \dots, K_d) \;=\; 
\frac 1{Z_{M,\bs i}} 
\prod_{j=1}^d  
\begin{pmatrix}
  |\Lambda_M| \\ K_j
\end{pmatrix}
\begin{pmatrix}
  |\Lambda_M| \\ K_j + I_j
\end{pmatrix}
\;.
\end{equation*}

For $0\le n\le |\Lambda_M|$ let
\begin{equation}
\label{g20}
h_n(a) = \Big ( \frac{|\Lambda_M| - a}{1+a}\Big) \Big(
\frac{|\Lambda_M| - (a+n)}{1+(a+n)}\Big )\; . 
\end{equation}
$h_n$ is a strictly convex, strictly decreasing function in the
interval $[0, |\Lambda_M| -n]$. Moreover, for $n<m$, $h_n(a) < h_m(a)$
on the interval $[0, |\Lambda_M| -m]$.

Denote by $\mf d_j$ the configuration of $\bb N^d$ with a unique
particle at coordinate $j$. An elementary computation shows that 
\begin{equation}
\label{g19}
\frac{\bar \nu (\bs K + \mf d_j)} {\bar \nu (\bs K + \mf d_k)}
\; \le 1\; \quad \text{if and only if} \quad
h_{I_j}(K_j) \le h_{I_k}(K_k)\;, 
\end{equation}
where summation is understood componentwise.

Denote by $\tilde {\bs K}$ an ordered solution of \eqref{g18} below and fix
a function $F$ in $L^2(\bar \nu_{\Lambda_M,\bs i})$. We have that
\begin{equation}
\label{g23}
E_{\bar \nu_{\Lambda_M,\bs i}}[F;F] \;\le\; 
E_{\bar \nu_{\Lambda_M,\bs i}} \big [ \{F(\bs K) - F(\tilde {\bs
  K})\}^2 \big ]\;.
\end{equation}

For each $\bs K$ in the hyperplane $\mc H$, consider the following
infinite path. Let $\bs K^0=\bs K$ and assume that $\bs K^0, \cdots,
\bs K^\ell$ have been defined. Let $j_0$, $k_0$ such that
\begin{equation}
\label{g22}
h_{I_{k_0}}(K^\ell_{k_0}-1) = \min_k h_{I_{k}}(K^\ell_{k}-1)\;, \quad
h_{I_{j_0}}(K^\ell_{j_0}) = \max_j h_{I_{j}}(K^\ell_{j})\;.
\end{equation}
If $\bs K^\ell$ is a solution of \eqref{g18}
($h_{I_{k_0}}(K^\ell_{k_0}-1) \ge h_{I_{j_0}}(K^\ell_{j_0})$), let
$\bs K^{\ell +1} = \bs K^\ell$; otherwise, let $\bs K^{\ell +1} = \bs
K^\ell - \mf d_{k_0} + \mf d_{j_0}$. In this latter case, by
\eqref{g21}, $\bar \nu_{\Lambda_M,\bs i} (\bs K^{\ell}) < \bar
\nu_{\Lambda_M,\bs i} (\bs K^{\ell +1})$.  Since $\mc H$ is finite
and since $\bar \nu_{\Lambda_M,\bs i} (\bs K^\ell)$ decreases
whenever $\bs K^\ell$ is not a solution of \eqref{g18}, the path
reaches eventually a solution.  The path can therefore be written as
$(\bs K^0 , \dots , \bs K^{\ell_0}, \bs K^{\ell_0}, \dots)$, where
$\bs K^{\ell_0}$ solves \eqref{g18} and $\bar \nu_{\Lambda_M,\bs i}
(\bs K^m) < \bar \nu_{\Lambda_M,\bs i} (\bs K^n)$ for $0\le m < n
\le \ell_0$.

By the end of the proof of Lemma \ref{gs11} below, there is a path
from $\bs K^{\ell_0}$ to $\tilde {\bs K}$ of length less than or equal
to $d/2$, passing only by solutions of \eqref{g18} and such that all
configurations visited have the same probability. Juxtaposing the two
previous paths, we obtain the path $\Gamma (\bs K, \tilde {\bs K}) =
(\bs K = \bs K^0 , \dots , \bs K^{\ell_0}, \dots , \bs K^{\ell_{\bs
    K}} = \tilde {\bs K})$, where $\ell = \ell_{\bs K}$ stands for the
total length of the path. By construction, the probability of the
configurations visited is non decreasing.

We are now ready to estimate the right hand side of \eqref{g23}. By
Schwarz inequality,
\begin{eqnarray*}
\!\!\!\!\!\!\!\!\!\!\!\! &&
E_{\bar \nu_{\Lambda_M,\bs i}} \big [ \{F(\bs K) - 
F(\tilde {\bs K})\}^2 \big ] \\
\!\!\!\!\!\!\!\!\!\!\!\! && \quad
\le\; \sum_{\bs K} \bar \nu_{\Lambda_M,\bs i}(\bs K) \ell_K
\sum_{j=0}^{\ell_K-1} \{F(\bs K_{j+1}) - F(\bs K_j)\}^2
\end{eqnarray*}
Since we just need a polynomial bound on the spectral gap and since
this method can not provide a sharp estimate, we bound the length of a
path $\ell_K$ by the total number of configurations $|\mc H| \le
|\Lambda_M|^d$. On the other hand, since $\bar \nu_{\Lambda_M,\bs i}
(\bs K) \le \bar \nu_{\Lambda_M,\bs i} (\bs K_j)$, we may replace
the former by the latter. Finally, inverting the order of summations
and estimating the total number of configurations which contains in
its path to $\tilde {\bs K}$ a fixed couple $\bs K_j$, $\bs K_{j+1}$
by the total number of configurations, we get that the previous
expression is less than or equal to
\begin{equation*}
C_0 M^{2 d^2} \sum_{\bs K} \sum_{\bs L \sim K} 
\bar \nu_{\Lambda_M,\bs i}(\bs K) \{F(\bs L) - F(\bs K)\}^2
\end{equation*}
for some universal constant $C_0$.  In this formula, the second sum is
carried over all configurations $\bs L$ which can be obtained from
$\bs K$ by letting a particle jump from a site to another: $\bs L =
\bs K - \mf d_{j} + \mf d_{k}$ for some $j\not = k$. By the explicit
formula for the Dirichlet form of $F$ derived above, this expression
is less than or equal to
\begin{equation*}
C_0 M^{3d + 2 d^2} \< F, -\mc L_M^c F\>_{\bar \nu_{\Lambda_M,\bs i}}\; .
\end{equation*}
It remains to apply Lemma \ref{gs12} to conclude the proof of the
spectral gap. \qed
\smallskip

We conclude this section with a result used in the proof of
Proposition \ref{gs10}.

\begin{lemma}
\label{gs11}
Fix $\bs i$ in $\mf V_M$ such that $I_1 \le \cdots \le
I_d$, $|\Lambda_M|\bs i:=(I_0,\dots,I_d)$. The system of equations
\begin{equation}
\label{g18}
\left\{
\begin{array}{l}
{\displaystyle \vphantom{\Big\{}
h_{I_k}(K_k-1) \ge h_{I_j}(K_j) \quad\text{for all $1\le k, j \le d$}} \\
{\displaystyle
2\sum_{j} K_j = I_0 - \sum_j I_j } \\
{\displaystyle
0\le K_j \le |\Lambda_M|\;, \quad 1\le j\le d\;,}
\end{array}
\right.
\end{equation}
has a solution such that $K_d \le \cdots \le K_1$. Moreover, if $\bs
K$, $\bs L$ are two solutions of \eqref{g18}, then $|K_j-L_j|\le 1$
for all $1\le j\le d$ and $\bar \nu_{\Lambda_M,\bs i}(\bs K) =
\bar \nu_{\Lambda_M,\bs i}(\bs L)$. 
\end{lemma}

\begin{proof}
To prove the existence of a solution, recall from \eqref{g19} that
\begin{equation}
\label{g21}
\frac{\bar \nu (\bs K + \mf d_j - \mf d_k )} {\bar \nu (\bs K )}
\; \le 1\; \quad \text{if and only if} \quad
h_{I_j}(K_j) \le h_{I_k}(K_k-1)\;. 
\end{equation}
Consider a configuration $\bs K^*$ which maximizes the probability
$\bar \nu_{\Lambda_M,\bs i}$. The inequality on the left hand side
of the previous displayed formula is satisfied for all $j$, $k$. In
particular, $\bs K^*$ solves \eqref{g18}.

Fix a solution of \eqref{g18}.  We claim that $K_j \le K_i$ if $I_i <
I_j$.  Assume by contradiction that $K_i < K_j$. In this case
\begin{equation*}
h_{I_i}(K_i) \;\le\; h_{I_j}(K_j-1) \;<\; h_{I_i}(K_j-1) \;\le\;
h_{I_i}(K_i)\;, 
\end{equation*}
which is a contradiction. Here, the first inequality follows from the
first property in \eqref{g18} of $\bs K$, the second from the fact
that $h_{I_j} < h_{I_i}$ and the last from the relation $K_i \le
K_j-1$.

Suppose that $I_i = I_j$ for some $i<j$ and that $K_i<K_j$ for a
solution $\bs K$ of \eqref{g18}. Let $\tilde {\bs K}$ be such that
$\tilde {\bs K}_k = \bs K_k$ for $k\not = i$, $j$; $\tilde {\bs K}_i =
\bs K_j$, $\tilde {\bs K}_j = \bs K_i$. It is easy to check that
$\tilde {\bs K}$ is also a solution of \eqref{g18}. This observation
together with the estimate derived in the previous paragraph show that
there exists a solution $\bs K$ of \eqref{g18} with $K_d \le \cdots
\le K_1$.

Finally, let $\bs K$, $\bs L$ be two solutions of \eqref{g18}.
Suppose by contradiction that $L_j \le K_j -2$ for some $j$. Since
$\sum_j K_j = \sum_j L_j$, there exists $i$ such that $K_i<L_i$. In
particular, 
\begin{equation*}
h_{I_i}(L_i-1) \;\le\; h_{I_i}(K_i) \;\le\; h_{I_j}(K_j-1)\;<\;
h_{I_j}(L_j) \;\le\; h_{I_i}(L_i-1)\;.
\end{equation*}
The first and third inequalities follow from the fact that $h_{I_i}$,
$h_{I_j}$ are strictly decreasing functions and the relations
$K_i<L_i$, $L_j<K_j-1$; while the second and fourth inequalities
follow from the property of $\bs K$, $\bs L$. This proves the first
property of $\bs K$, $\bs L$.

To prove the second property of $\bs K$, $\bs L$, consider a path from
$\bs K$ to $\bs L$: $\bs K= \bs M^0, \dots, \bs M^\ell = \bs L$, for
each $0\le i <\ell$, $\bs M^{i+1} = \bs M^i +\mf d_j - \mf d_k$ for
some $j=j(i)$, $k=k(i)$. It is not difficult to show that there exists
such a path with $\ell \le d/2$, $\bs M^i$ solving \eqref{g18} for all
$i$.

Fix $i$ and let $\bs M^i = \bs M$, $\bs M^{i+1} = \bs M+\mf d_j - \mf
d_k$, $\bs M^i$, $\bs M^{i+1}$ solving \eqref{g18}. Since $\bs M$
solves \eqref{g18}, $h_{I_j} (M_j) \le h_{I_k}(M_k-1)$.  Using now
that $\bs M+\mf d_j - \mf d_k$ solves \eqref{g18}, we obtain the
reverse inequality so that $h_{I_j} (M_j) = h_{I_k}(M_k-1)$. In
particular, in view of \eqref{g21}, $\bar \nu_{\Lambda_M,\bs i}(\bs
M^i) = \bar \nu_{\Lambda_M,\bs i}(\bs M^{i+1})$.  This concludes
the proof of the lemma.
\end{proof}

\section{Equivalence of ensembles}
\label{equiv}

We prove in this section the equivalence of ensembles for the
stochastic lattice gas introduced in Section \ref{sec1}. Recall the
definition of the set $\mf V_L$ and of the canonical measures
$\nu_{\Lambda_L, \bs i}$. Notice that for every $\bs \lambda$ in $\bb
R^{d+1}$
\begin{equation*}
\nu_{\Lambda_L,\bs i} (\cdot ) \;=\; \mu^{\Lambda_L}_{\bs \lambda} 
\big( \, \cdot \, \big\vert \, \bs I^L =\bs i \big)\;.
\end{equation*}

For $(\rho,\bs p)$ in $\mf A$ the expectation of the one site random
variable $\bs I (\eta_x)$ under the product measure
$\mu^{\Lambda_L}_{\Lambda(\rho,\bs p)}$ is equal to $(\rho,\bs p)$. It
defines a map from $\mf A$ to the set of probability measures on
$\big(\{0,1\}^{\mc V}\big)^{\Lambda_L}$. Since this map is uniformly
continuous it may be extended continuously to the closure of $\mf A$.
For each $\bs i \in\mf V_L$ denote by $\mu^L_{\Lambda(\bs i)}$ the
corresponding product measure by this map. Hence we have a one to one
correspondence between the canonical measures $\{\nu_{\Lambda_L,\bs i}
: \bs i\in\mf V_L \}$ and the so-called grand canonical measures
$\{\mu^N_{\bs \Lambda(\bs i)} : \bs i\in\mf V_L\}$.  

Let $\< g; f\>_{\mu}$ stand for the covariance of $g$, $f$ with
respect to $\mu$: $\< g; f\>_{\mu} = E_{\mu}[fg] - E_{\mu}[f]
E_{\mu}[g]$ and $\< f,g\>_{\mu}$ for the inner product in $L^2(\mu)$.

\begin{proposition}
\label{hs1}
Fix a cube $\Lambda_\ell\subset\Lambda_L$. For each $\bs i\in \mf V_L$
denote by $\nu^{\ell}$ the projection of the canonical measure
$\nu_{\Lambda_L,\bs i}$ on $\Lambda_\ell$ and by $\mu^{\ell}$ the
projection of the grand canonical measure $\mu^L_{\bs \Lambda (\bs
  i)}$ on $\Lambda_\ell$. Then, there exists a finite constant
$C(\ell, \mc V)$, depending only on $\ell$ and $\mc V$, such that
\begin{equation*}
\big|\,E_{\mu^{\ell}}[f]-E_{\nu^{\ell}}[f]\,\big|
\le \frac{C(\ell, \mc V)}{|\Lambda_L|}\;\< f; f\>^{1/2}_{\mu^{\ell}}
\end{equation*}
for every $f:\big(\{0,1\}^{\mc V}\big)^{\Lambda_{\ell}}\mapsto \bb R$.
\end{proposition}

\begin{proof}
Since $\nu^{\ell}$ is absolutely continuous with respect to
$\mu^{\ell}$, by Schwarz inequality, 
\begin{eqnarray*}
\big|\,E_{\nu^{\ell}}[f]-E_{\mu^{\ell}}[f]\,\big|& = &
\Big|E_{\mu^{\ell}}\big[\, \big( \frac{d\nu^{\ell}}{d\mu^{\ell}}\big)
\big( f - E_{\mu^{\ell}}[f] \big) \, \big]\Big| \nonumber\\ 
& \le & \< \;\frac{d\nu^{\ell}}{d\mu^{\ell}}\; ;\;
\frac{d\nu^{\ell}}{d\mu^{\ell}} \; \>^{1/2}_{\mu^{\ell}} \; \< f;
f\>^{1/2}_{\mu^{\ell}}.
\end{eqnarray*}
Since $\mu^L_{\bs \Lambda(\bs i)}$ is a product measure, for any $\xi$
in $\big(\{0,1\}^{\mc V}\big)^{\Lambda_{\ell}}$,
\begin{equation*}
\frac{d\nu^{\ell}}{d\mu^{\ell}}(\xi) = \frac{ \mu^L_{\bs \Lambda (\bs
    i)}\big[ \;\sum_{x\in\Lambda_L \setminus \Lambda_{\ell}}\bs
  I(\eta_x) = |\Lambda_L|\bs i - |\Lambda_\ell| \bs I^{\ell}(\xi) \;\big] }{
  \mu^L_{\bs \Lambda (\bs i)}\big[\; \sum_{x\in\Lambda_L}\bs
  I(\eta_x) = |\Lambda_L| \bs i \;\big] } \;,
\end{equation*}
where $\bs I^{\ell}(\xi) = |\Lambda_\ell|^{-1} \sum_{x\in\Lambda_\ell}
\bs I(\xi_x)$.  Under $\mu^L_{\bs \Lambda (\bs i)}$, $\bs I(\eta_x)$
are i.i.d. random variables taking a finite number of values. By
Theorem VII.12 in \cite{p}, there exists a finite constant $C_0(\ell,
\mc V)$, depending only on $\ell$ and $\mc V$, such that,
\begin{equation*}
\big| \frac{d\nu^{\ell}}{d\mu^{\ell}}(\xi) -1 \big| \le 
C_0(\ell, \mc V) \frac 1{|\Lambda_L|}
\end{equation*}
uniformly in $\xi$. This concludes the proof of the lemma.
\end{proof}

\medskip
\noindent{\bf Acknowledgment} The authors would like to thank
H. T. Yau for very fruitful discussions and his warm hospitality at
Stanford University.

\end{document}